\documentclass[12pt]{thesis}  
\usepackage{graphicx}
\usepackage{color}
\usepackage{transparent}
\usepackage{lscape}
\usepackage{indentfirst}
\usepackage{latexsym}
\usepackage{multirow}
\usepackage{tabls}
\usepackage{wrapfig}
\usepackage{slashbox}
\usepackage{longtable}
\usepackage{enumerate}
\usepackage{tikz}
\usetikzlibrary{matrix,arrows}
\renewcommand{\baselinestretch}{2}
\newcommand{\wt}{\widetilde}
\newcommand{\mc}{\mathcal{C}}

\newcommand{\mg}{\mathcal{G}}
\newcommand{\mb}{\mathcal{B}}
\newcommand{\bdyinf}{\partial_\infty}
\newcommand{\rc}{\mathbf{R}^C}
\newcommand{\ep}{\epsilon}

\newtheorem{mythm}{Theorem}
\newtheorem{myprop}{Proposition}[chapter]
\newtheorem{mylemma}[myprop]{Lemma}
\newtheorem{mycor}[myprop]{Corollary}

\newcommand{\qed}{\nobreak \ifvmode \relax \else
      \ifdim\lastskip<1.5em \hskip-\lastskip
      \hskip1.5em plus0em minus0.5em \fi \nobreak
      \vrule height0.75em width0.5em depth0.25em\fi} 

\begin{document}


\hbox{\ }

\renewcommand{\baselinestretch}{1}
\small \normalsize

\begin{center}
\large{{ABSTRACT}}

\vspace{3em}

\end{center}
\hspace{-.15in}
\begin{tabular}{ll}
Title of dissertation:    & {\large  LENGTH SPECTRAL RIGIDITY }\\
&				      {\large  OF NON-POSITIVELY CURVED SURFACES} \\
\ \\
&                          {\large  Jeffrey Frazier, Doctor of Philosophy, 2012} \\
\ \\
Dissertation directed by: & {\large  Professor William Goldman} \\
&  				{\large	 Department of Mathematics } \\
\end{tabular}

\vspace{3em}

\renewcommand{\baselinestretch}{2}
\large \normalsize

Length spectral rigidity is the question of under what circumstances the geometry of a surface can be determined, up to isotopy, by knowing only the lengths of its closed geodesics. It is known that this can be done for negatively curved Riemannian surfaces, as well as for negatively-curved cone surfaces. Steps are taken toward showing that this holds also for flat cone surfaces, and it is shown that the lengths of closed geodesics are also enough to determine which of these three categories a geometric surface falls into. Techniques of Gromov, Bonahon, and Otal are explained and adapted, such as topological conjugacy, geodesic currents, Liouville measures, and the average angle between two geometric surfaces. 

\thispagestyle{empty}
\hbox{\ }
\vspace{1in}
\renewcommand{\baselinestretch}{1}
\small\normalsize
\begin{center}

\large{{LENGTH SPECTRAL RIGIDITY OF \\
NON-POSITIVELY CURVED SURFACES}}\\
\ \\
\ \\
\large{by} \\
\ \\
\large{Jeffrey Frazier}
\ \\
\ \\
\ \\
\ \\
\normalsize
Dissertation submitted to the Faculty of the Graduate School of the \\
University of Maryland, College Park in partial fulfillment \\
of the requirements for the degree of \\
Doctor of Philosophy \\
2012
\end{center}

\vspace{7.5em}

\noindent Advisory Committee: \\
Dr. William Goldman, Chair/Advisor \\
Dr. Scott Wolpert \\
Dr. Karin Melnick \\
Dr. Richard Wentworth \\
Dr. Dieter Brill \\

\thispagestyle{empty}
\hbox{\ }

\vfill
\renewcommand{\baselinestretch}{1}
\small\normalsize

\vspace{-.65in}

\begin{center}
\large{\copyright \hbox{ }Copyright by\\
Jeffrey Frazier  
\\
2012}
\end{center}

\vfill

\pagestyle{plain}
\pagenumbering{roman}
\setcounter{page}{2}

\renewcommand{\baselinestretch}{1}
\small\normalsize
\tableofcontents 
\newpage
\listoffigures 
\newpage

\newpage
\setlength{\parskip}{0em}
\renewcommand{\baselinestretch}{2}
\small\normalsize

\setcounter{page}{1}
\pagenumbering{arabic}

\renewcommand{\thechapter}{1}

\chapter{Introduction}

Let $S$ be a fixed closed, orientable surface of genus at least 2. Let $\Gamma$ be the fundamental group of $S$, and $C$ the set of isotopy classes of closed curves on $S$.

Consider the following three moduli spaces of marked geometric surfaces homeomorphic to $S$, each up to isometry isotopic to the identity:

\begin{itemize}
\item Neg($S$) - Riemannian surfaces of variable, but strictly negative curvature,
\item Neg*($S$) - Riemannian cone surfaces of strictly negative curvature, with all cone angles in excess of $2\pi$,
\item Flat*($S$) - Flat cone surfaces, with all cone angles in excess of $2\pi$.
\end{itemize}

Let NonPos($S$) denote the disjoint union of these three moduli spaces. Any $X \in \textrm{NonPos}(S)$ defines a {\em marked length spectrum}, a functional on $C$ which associates to each $\alpha \in C$ the length of the unique geodesic in $X$ which is in the isotopy class $\alpha$. Letting $\mathbf{R}^C$ denote the space of all functionals on $C$, this determines a mapping
$$l: \textrm{NonPos}(S)\longrightarrow \mathbf{R}^C.$$

The question of spectral rigidity asks whether this mapping, possibly when restricted to some proper subset, is injective. If $l$ is injective on some subspace $\Sigma \subset \textrm{NonPos}(S)$, then $\Sigma$ is said to be {\em spectrally rigid} over $C$. This means that knowing the lengths of all closed geodesics on a surface in $\Sigma$ is enough to determine the entire geometry of the surface up to isotopy.

This paper proves the following two results, with an immediate corollary:

\vspace{10mm}

\noindent\textbf{Theorem 1} The image of Neg($S$) in $\mathbf{R}^C$ intersects neither the image of Neg*($S$) nor the image of Flat*($S$). That is, no Riemannian surface has the same length spectrum as a surface with cone points.

\vspace{10mm}

\noindent\textbf{Theorem 2} The images of Neg*($S$) and Flat*($S$) in $\mathbf{R}^C$ are disjoint. That is, no negatively curved cone surface has the same length spectrum as a flat surface.

\vspace{10mm}

\noindent\textbf{Corollary} These three moduli spaces have pairwise disjoint images in $\rc$.

\vspace{10mm}

J.P. Otal \cite{Otal} proved in 1990 that Neg($S$) is spectrally rigid, i.e. that $l$ is injective when restricted to Neg($S$). Sa'ar Hersonsky and Fr\'ed\'eric Paulin \cite{Hersonsky} adapted Otal's proof in 1997 to show that Neg*($S$) is also spectrally rigid. If it were proved that Flat($S$) were spectrally rigid as well, then these results, combined with Theorems 1 and 2, would show that $l$ is injective on all of NonPos($S$). To the author's knowledge, this has not been proved. An discussion is included as to why the author believes this result should be true, as well as a framework for a possible argument, but no complete proof is given.

A recent paper of Moon Duchin, Christopher Leininger, and Kasra Rafi \cite{Duchin} showed that the subspace of Flat*($S$) consisting of all flat surfaces whose cone angles are multiples of $\pi$ is spectrally rigid over all {\em simple} closed curves, thus proving a much stronger rigidity result for a much smaller class of structures. Their proof relies strongly on the fact that such a flat surface can be defined from a quadratic differential on a Riemann surface, which is not true of a general flat surface.

The method of proof of Theorems 1 and 2 uses machinery developed by Bonahon (\cite{Bonahon1}, \cite{Bonahon2}) and Otal to translate between length spectra and measure-theoretic objects called {\em Liouville measures}, which are definable from a surface (see prop. \ref{spectrumequiv}, and also see \cite{Ledrappier} for other types of objects which are equivalent to length spectra). A large role is also played by the notion of topological conjugacy, by which the geodesic structures of any two surfaces in NonPos($S$) may be identified with each other (prop. \ref{bigconjugacy}).

The paper is structured as follows. Chapter 2 outlines the theory of conjugacy and Liouville currents for Neg($S$). Chapters 3 and 4 do the same for Neg*($S$) and Flat*($S$), focusing on the differences that arise from the addition of cone points and flat curvature. Chapter 5 presents the proofs of Theorems 1 and 2, as well as arguments concerning the rigidity of Flat*($S$). 

\renewcommand{\thechapter}{2}

\chapter{Negatively curved Riemannian surfaces}

Let $X$ be a negatively curved Riemannian surface with a marking homeomorphism $X\rightarrow S$. Throughout, the Riemannian surface $X$ and the marking may vary, but $S$ will be forever fixed.

A result of Melvyn Berger \cite{Berger} from the early 70's gives a complete classification of all compact, negatively curved Riemannian surfaces in terms of three pieces of data: the genus of the surface, the conformal class of the metric, and the pointwise Gaussian curvature function.

\begin{myprop}
\label{berger}
Let $\Sigma$ be a compact Riemann surface. Then any smooth, negative function $K:\Sigma\rightarrow \mathbf{R}$ is the Gaussian curvature of a unique metric in the conformal class of $\Sigma$.
\end{myprop}

The proof involves choosing a base metric in the conformal class and solving an elliptic partial differential equation to find a function which conformally deforms the base metric to one with the given curvature. Similar results for cone structures will be noted in further chapters.

With this result, the moduli space of negatively curved Riemannian structures on $S$ can be described as a product of the Teichmuller space of conformal structures with a function space of smooth, negative functions on $S$.

\section{Boundary at infinity and conjugacy}

The notion of topological conjugacy of surfaces is central to the arguments herein. The idea is that given any two homeomorphic non-positively curved Riemannian surfaces, there is a direct correspondence between the structures of their geodesics. This correspondence will be described below for negatively curved Riemannian surfaces, and extended to cone surfaces in the following chapters. First we must define the {\em boundary at infinity} of a negatively curved Riemannian surface $X$.

\label{boundaryatinfinity}
Let $\wt{X}$ denote the metric universal cover of $X$. Two oriented geodesics of $\wt{X}$ are said to be {\em asymptotic} if they stay within a bounded distance of each other for all positive time (note that this is independent of orientation-preserving reparametrization). This is clearly an equivalence relation.

The space of equivalence classes, with a cone topology defined from half-spaces in $\wt{X}$ (see \cite{Bridson}), is called the boundary at infinity of $\wt{X}$ and is denoted $\bdyinf \wt{X}$, or by slight abuse, $\bdyinf X$. A single equivalence class in $\bdyinf X$ is a ``point" at infinity. Topologically $\bdyinf X$ is a circle, and the union $\wt{X}\cup\bdyinf X$ can be topologized so that it is a compactification of $\wt{X}$, homeomorphic to a 2-disk. Since $\pi_1(X)$ acts on $\wt{X}$ by isometries, this action extends to $\bdyinf X$.

An important feature of a negatively curved Riemannian surface $X$ is that $\bdyinf X$ is homeomorphic to the ``visual sphere" at any point in $\wt{X}$, via the exponential map. If $T^1_p \wt{X}$ is the circle of unit vectors at a point $p \in \wt{X}$, then the map $T^1_p \wt{X} \rightarrow \bdyinf X$ which sends a vector $v$ to the asymptotic class of the geodesic through $p$ in the direction of $v$ is a homeomorphism. This will contrast with the cone structures considered later.

\begin{figure}[htb]
  \centering
  \def\svgwidth{150pt}
  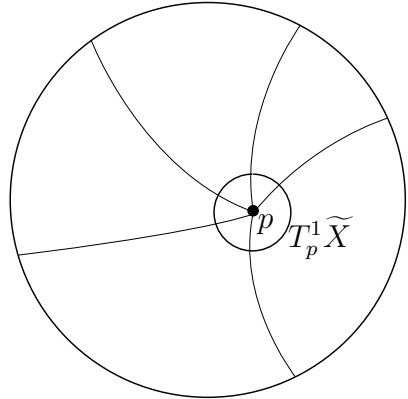
  \caption{Exponential map is a homeomorphism to the boundary}
\end{figure}

Any oriented geodesic $\gamma$ on $\wt{X}$ determines two endpoints at infinity, one in the positive direction and one in the negative. Conversely, any two distinct points $a$ and $b$ at infinity determine a unique oriented geodesic with initial point $a$ and terminal point $b$ (note that this latter property will fail for the flat structures considered later).

The universal cover $\wt{X}$ is an example of a {\em Hadamard space}: a complete, simply connected, non-positively curved length space. It is also a hyperbolic space in the sense of Gromov (see \cite{Gromov}). The group $\Gamma = \pi_1(S)$ acts on $\wt{X}$ via the marking $X\rightarrow S$. This action is discrete and cocompact, as a group of isometries. For the following result, see \cite{Gromov} and \cite{Croke}.

\begin{myprop}
\label{conjugacy}
If a finitely generated group $G$ acts discretely, cocompactly, and isometrically on two hyperbolic Hadamard spaces, then there is a $G$-equivariant homeomorphism between their respective boundaries at infinity.
\end{myprop}

\begin{mycor}
Let $X_1$ and $X_2$ be negatively curved Riemannian surfaces with markings $f_i:X_i\rightarrow S$. Lift $f_2^{-1}\circ f_1$ to the universal covers and extend to the boundary to obtain a map $\phi: \bdyinf X_1 \rightarrow \bdyinf X_2$. Then $\phi$ is a homeomorphism equivariant with respect to the actions of $\Gamma$.
\end{mycor}

The map $\phi$ in the proposition is called a {\em conjugacy map}. Note that, since $S$ is a closed surface, the endpoints of invariant axes of isometries in $\Gamma$ are dense in the boundary. This implies that $\phi$ must be unique, since the requirement that $\phi$ is $\Gamma$-equivariant determines the map on this dense subset.

Also note that the conjugacy map depends on the markings chosen for $X_1$ and $X_2$, since these markings are needed to define the actions of $\Gamma = \pi_1(S)$ on the boundaries. If one or both markings are changed, the conjugacy map will be altered by an appropriate mapping class. In fact, the conjugacy map only depends on the isotopy types of the markings, since two isotopic maps $X_1\rightarrow X_2$ will have the same extension to the boundary. Thus the natural setting is the space Neg($S$) of isotopy classes of negatively curved manifolds marked to $S$.

The importance of the conjugacy map is that (as we will see) there are various objects definable from a negatively curved surface using only the action of $\Gamma$ on its boundary. Then the existence of the conjugacy map will show that such objects depend (up to topological considerations) only on the underlying smooth surface.

Conjugacy maps can be better understood by thinking about geodesics and quasi-geodesics on a universal cover $\wt{X}$ (a quasi-geodesic is the image of quasi-isometry $\mathbf{R}\rightarrow \wt{X}$). It is a well-known fact that any quasi-geodesic in a simply-connected, negatively curved space is a bounded distance from a unique geodesic (this will contrast slightly with the flat case). The composition $f = f_2^{-1}\circ f_1$ of the markings is a quasi-isometry of $X_1$ and $X_2$, simply because they are compact. Thus it lifts to a quasi-isometry $\wt{f}:\wt{X_1}\rightarrow\wt{X_2}$.

Then the conjugacy map between $\bdyinf X_1$ and $\bdyinf X_2$ can be defined as follows. Choose an oriented geodesic $\gamma$ on $\wt{X_1}$ and let $a$ be its terminal endpoint at infinity. The ``same" curve in $\wt{X_2}$ (in other words, $\wt{f}\circ\gamma$) is probably not a geodesic in the metric of $\wt{X_2}$, but it is a quasi-geodesic, and its terminal endpoint at infinity is defined to be $\phi(a)$. It is clear that a mapping defined this way is $\Gamma$-equivariant.

It will also be useful to have a purely topological model for the boundary at infinity, defined using only the smooth surface $S$. Let $G$ be the Cayley graph of $\Gamma$, based on the standard set of generators. Define $\bdyinf S$ to be the set of asymptotic classes of graph geodesics in $G$ (note that two graph geodesics are asymptotic if and only if they are eventually equal). As before, the action of $\Gamma$ on $G$ extends to $\bdyinf S$, since this action is by graph isometries.

Given a marking of a Riemannian manifold $X_1\rightarrow S$, $G$ can be embedded into $\wt{X_1}$ after choosing a basepoint on the surface, and such an embedding is a quasi-isometry from the graph metric on $G$ to the Riemannian metric on $\wt{X_1}$. Then, as above, any geodesic in $G$ becomes a quasi-geodesic in $\wt{X_1}$, so by similar considerations there is a $\Gamma$-equivariant map $\phi_1: \bdyinf X_1\rightarrow\bdyinf S$. Given another structure $X_2$, there is a similar map $\phi_2: \bdyinf X_2 \rightarrow \bdyinf S$, and $\phi = \phi_2^{-1}\circ\phi_1$ is the conjugacy map between the two structures.

\section{Geodesic currents}

The space of geodesic currents was introduced by Bonahon (\cite{Bonahon1}, \cite{Bonahon2}) as a completion of the set of weighted closed curves on a surface. Bonahon showed that the Fricke space of hyperbolic metrics on $S$ can be embedded naturally into the space of geodesic currents; this has been extended to more general types of moduli spaces of geometric structures, via the Liouville measures defined below.

Given a negatively curved Riemannian surface $X$, we can identify the space $\mathcal{G}(\wt{X})$ of complete unoriented geodesics of $\wt{X}$ with $(\bdyinf X \times \bdyinf X \setminus \Delta)/\mathbf{Z}_2$; that is, unordered pairs of distinct points at infinity. Since $\bdyinf X$ is topologically a circle, $\mathcal{G}(\wt{X})$ is homeomorphic to an open Moebius strip. A metric geodesic current on $X$ is defined to be a Borel measure on $\mathcal{G}(\wt{X})$ which is invariant under the (diagonal) action of $\pi_1(X)$. The collection of all such measures is denoted $\mc(X)$, and given a weak* uniform measure topology.

This construction can be mimicked using only the topological surface $S$. Define $\mg(\wt{S})$ as $(\bdyinf S \times \bdyinf S \setminus \Delta)/\mathbf{Z}_2$, and let $\mc(S)$ denote the space of topological geodesic currents on $S$, i.e. $\Gamma$-invariant Borel measures on $\mg(\wt{S})$. Given a marking $f:X\rightarrow S$, the conjugacy map $\phi:\bdyinf X\rightarrow\bdyinf S$ induces an equivariant homeomorphism $\phi\times\phi:\mg(\wt{X})\rightarrow\mg(\wt{S})$. Thus $\mc(X)$ is identified with $\mc(S)$. Note once again that this identification depends on the isotopy type of the marking.

A geodesic current on $X$ can also be thought of as a transverse invariant measure to the geodesic foliation on the unit tangent bundle of $X$. In fact, this is how Bonahon first defined currents in \cite{Bonahon1}. The geodesic foliation has codimension 2, so the invariant measures are defined on subsurfaces. Since this foliation can be constructed using only the $\Gamma$-action on $\bdyinf X$, the currents so defined again depend only on $S$ and the marking.

Recall that Neg($S$) denotes the moduli space of isotopy classes of negatively curved manifolds with markings to $S$. If $X$ denotes a class in Neg($S$), then we can write $\mc(X)$ without any ambiguity, since for any two marked surfaces in the class $X$, the conjugacy between their boundaries is induced by an isometry. Furthermore, $\mc(X)$ is {\em uniquely} identified with $\mc(S)$, since the identification depends only on the isotopy type of the marking. This means we can pass freely and without comment between the (metric) currents $\mc(X)$ and the (topological) currents $\mc(S)$. We will take this point of view hereafter.

Let $C$ be the (discrete) set of homotopy classes of closed curves on $S$. For any class $\alpha \in C$, we can define a geodesic current (also denoted $\alpha\in\mc(S)$) as follows. Choose any $X \in \textrm{Neg}(S)$ , and let $\gamma$ be the unique $X$-geodesic in the class of $\alpha$. The complete lift of $\gamma$ to $\wt{X}$ can be thought of as a discrete subset of $\mathcal{G}(X)$, so define the current associated to $\alpha$ as the Dirac measure on this subset. This measure is trivially $\Gamma$-invariant, since the complete lift was taken. It is clear that the current thus defined does not depend on the choice of $X$.

\label{curvesdenseincurrents}
Extend this to the space $C\times\mathbf{R_+}$ of weighted curves by multiplying the Dirac measure by the weight. This gives an embedding of $C\times\mathbf{R_+}$ into $\mc(S)$, and one of the fundamental results about currents is that the image of this embedding is dense in $\mc(S)$ \cite{BonahonThesis}.

The space $\mc(S)$ also carries a symmetric, continuous bilinear form called the intersection form, which is an extension of the geometric intersection number on $C$. We will only need to intersect currents in the case that at least one current comes from $C$, so we give the definition in this case. Let $\alpha \in C$, and $\mu\in\mathcal{C}(S)$. Recall that $\alpha$ can be thought of as a conjugacy class of $\Gamma$, and let $\gamma\in\alpha$ be a representative. Choose a reference metric $X\in\textrm{Neg}(S)$, and let $I$ be a fundamental domain for the action of $\gamma$ on its axis in $\wt{X}$. Then $i(\mu,\alpha)$ is the $\mu$-measure of the set of all geodesics in $\mathcal{G}(\wt{X})$ which intersect $I$. Since $\mu$ is $\Gamma$-invariant, it is easy to see that this does not depend on the choice of conjugacy class representative or fundamental domain. It is also easy to see that the intersection number does not depend on the negatively curved reference metric chosen.

The following result of Otal \cite{Otal} will be important for much of what follows.

\begin{myprop}
\label{otalcurrents}
A geodesic current is determined by its intersection numbers with all currents in $C$. That is, if $\mu, \nu \in \mathcal{C}(S)$ and $i(\mu,\alpha) = i(\nu, \alpha)$ for all $\alpha \in C$, then $\mu = \nu$.
\end{myprop}

This result can be interpreted as follows. Let $\mathbf{R}^C$ denote the space of real-valued functionals on $C$. There is a mapping $I:\mc(S)\rightarrow \mathbf{R}^C$ given by $\mu \mapsto (\alpha \mapsto i(\mu, \alpha))$. The above proposition states that $I$ is injective.

\section{The Liouville current of a negatively curved metric}

Given any negatively curved Riemannian surface with a marking to $S$, there is an associated geodesic current in $\mc(S)$, called the Liouville current of the marked surface. The Liouville current of a marked surface depends on the marking up to isotopy, so this can be thought of as a mapping $L:\textrm{Neg}(S)\rightarrow \mc(S)$, denoted as $L(X) = L_X$. The Liouville current can be constructed in three ways, each of which is illuminating in different contexts. An outline of an argument is then given as to why each of these constructions yields the same result.

\subsection{The Liouville current as a transverse measure}

Recall that a geodesic current can be thought of as a transverse invariant measure to the geodesic foliation on the unit tangent bundle, with respect to some reference metric class. The topological structure of this foliation does not depend on the negatively curved reference metric chosen, so any convenient choice may be made.

Given a negatively curved Riemannian surface $X$, there is a standard volume measure on $T_1 X$ which is locally the product of the Riemannian volume form on $X$ with the Lebesgue measure on each circular fiber, and which is invariant under the geodesic flow. Taking the interior product of this 3-form with the unit vector field that generates the geodesic flow results in a 2-form which is zero along the leaves of the geodesic foliation. Then the absolute value of this 2-form is a transverse invariant measure, and this is the Liouville current $L_X$.

\subsection{The Liouville current from a cross-ratio on $\bdyinf S$}
\label{crossratiodef}
The Liouville current can also be constructed as a measure on $\mathcal{G}(\wt{S}) = (\bdyinf S\times\bdyinf S \setminus \Delta)/\mathbf{Z}_2$ by specifying the measure of each product rectangle. Let $[a,b]$ and $[c,d]$ be non-overlapping segments in $\bdyinf S$. Then $[a,b]\times[c,d]$ is a rectangle in $\mathcal{G}(\wt{S})$, consisting of all geodesics with one endpoint on $[a,b]$ and the other on $[c,d]$. Then we define $L_X([a,b]\times[c,d]) = |C_X(a,b,c,d)|$, where $C_X$ is a real-valued function called the {\em cross-ratio} of $X$ that takes four distinct points on $\bdyinf S$. It remains to define this cross-ratio.

\begin{figure}[htb]
  \centering
  \def\svgwidth{200pt}
  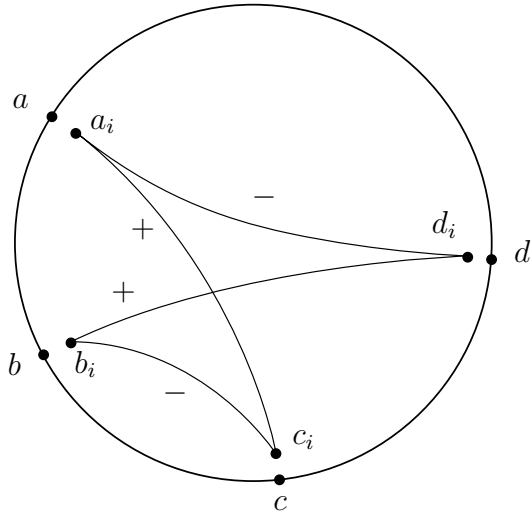
  \caption{The cross-ratio $C_X(a,b,c,d)$.}
\end{figure}

Let ${a_i}$ be a sequence of points in $\wt{X}$ limiting to $a$ at the boundary, and similarly for the other points $b,c,d$. Then the cross-ratio of $X$ is defined as \label{crossratio}
$$C_X(a,b,c,d) = \frac{1}{2}\lim (d(a_i, c_i) + d(b_i, d_i) - d(a_i, d_i) - d(b_i, c_i)),$$

\noindent where $d$ is the Riemannian distance function in $\wt{X}$, and the limit is taken as $a_i \rightarrow a$, etc. It can be shown using horocycles that this limit always exists and is finite, and does not depend on the choice of sequences limiting to $a,b,c,d$. This is done explicitly in section \ref{flatcrossratio} for the case of a flat structure. Note that $C_X$ is invariant under the action of $\Gamma$ on $\bdyinf S$, since $\Gamma$ acts by isometries.

Now that the measure $L_X$ is defined on rectangles, it can be extended to a Borel measure using the Carath\'eodory construction. This is carried out in detail in (\cite{Hersonsky} Thm 4.4). The basic steps are to use the rectangle measure to build an outer measure on $\mathcal{G}(\wt{S})$, invoke the Carath\'eodory construction to produce a collection of measurable sets, and then check that all Borel sets are measurable and that the measure produced does what we expect on the rectangles.

The resulting Borel measure on $\mathcal{G}(\wt{S})$ is invariant under $\Gamma$ since $C_X$ is, and this is the Liouville current $L_X$.

\subsection{The Liouville current in geodesic-angle coordinates}

Let $\gamma$ be a complete geodesic in $\wt{X}$, with endpoints $a, b\in \bdyinf S$, and let $\mg(\gamma)$ be the open set of all geodesics in $\mathcal{G}(\wt{S})$ that intersect $\gamma$ transversely in $\wt{X}$. Let $t\mapsto \gamma(t)$ be a unit-speed parametrization of $\gamma$.

\begin{figure}[htb]
  \centering
  \def\svgwidth{400pt}
  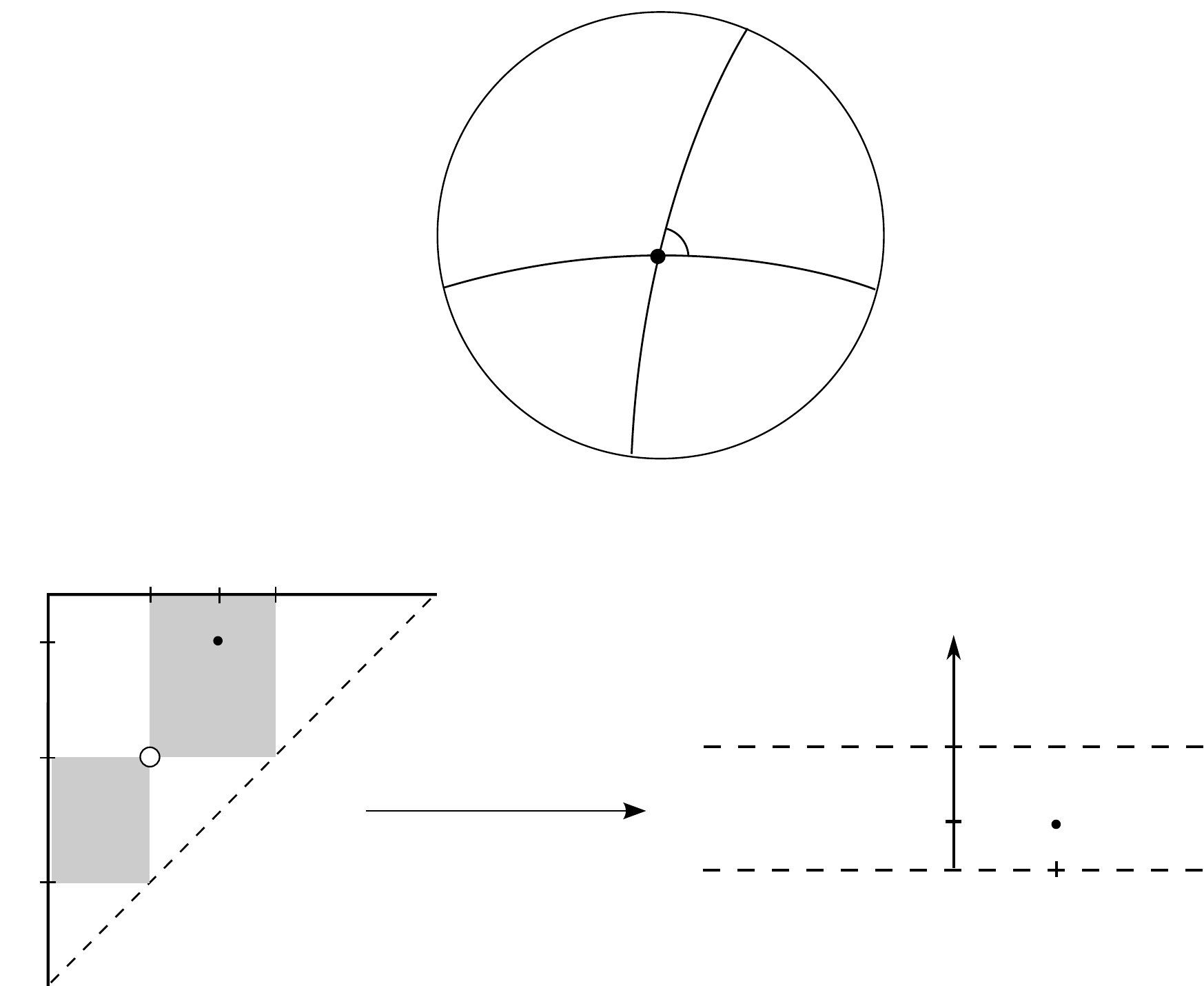
  \caption{A geodesic-angle coordinate on $\mg(\wt{S})$.}
\end{figure}

Any geodesic in $\mg(\gamma)$ intersects $\gamma$ at a single point, and with some angle. Then the mapping $\xi_{X,\gamma}:\mg(\gamma) \rightarrow \mathbf{R}\times(0,\pi)$, which takes a geodesic to $(t,\theta)$, the parameter of its intersection point with $\gamma$ together with the angle of intersection, is a homeomorphism. This map is called a {\em geodesic-angle coordinate} on $\mathcal{G}(\wt{S})$, and is uniquely defined by choosing the geodesic $\gamma$, along with an orientation and an origin for $\gamma$.

As $\gamma$ varies over all geodesics in $\wt{X}$, the sets $\mg(\gamma)$ form an open cover of $\mg(\wt{S})$. Each set in the open cover comes with a homeomorphism to $\mathbf{R}\times(0,\pi)$. Let $d\lambda = \frac{1}{2}\sin\theta d\theta dt$ be a measure on $\mathbf{R}\times(0,\pi)$, and pull back $d\lambda$ through each $\xi_{X,\gamma}$. This defines a measure on each $\mg(\gamma)$; to see that these measures agree on intersections, see (\cite{Santalo}, Ch. 19). Note that since $d\lambda$ is invariant in $t$, the choice of orientation and origin of $\gamma$ does not affect the pullback of the measure.

This defines a measure on $\mg(\wt{S})$, and again the fact that $\Gamma$ acts by isometries implies that this measure is $\Gamma$-invariant, so it is a geodesic current.

\subsection{All three constructions yield the same geodesic current $L_X$}
\label{constructionsequal}
It is important to note where, in each of these constructions, the specific geometry of the structure $X$ comes into play.

In the invariant transverse measure construction, it is in the use of the volume form on $X$, which is computed from the metric tensor. The derivative of the geodesic flow also depends strongly on the geometry of $X$, even though the foliation structure does not.

In the cross-ratio construction, it is in the cross-ratio function itself, which uses the Riemannian distance function on $\wt{X}$.

In the geodesic-angle construction, all the geometry of $X$ is contained in the coordinate functions $\xi_{X,\gamma}$. The measure $d\lambda$ can be thought of as the ``raw material" that all Liouville currents are made of, and the geodesic-angle coordinates describe how to arrange the raw materials to create the current that represents the specific Riemannian surface.

It remains to show that all three of these constructions result in the same geodesic current. The basis for this is the following lemma, the proof of which is outlined with references.

\begin{mylemma}
\label{lengthismeasure}
For each of the three constructions, the current $L_X$ produced satisfies the following property: given any geodesic segment $I$ in $\wt{X}$, the $L_X$-measure of the open set of all geodesics intersecting $I$ is equal to the length of $I$.
\end{mylemma}
{\em Proof.} This is easiest for the geodesic-angle construction, where it follows from a simple integration on $\mg(\gamma)$, where $\gamma$ is the complete geodesic carrying $I$, and using the fact that the parametrization of $\gamma$ is by arc-length.

For the invariant transverse measure construction, the argument can be found in \cite{Bonahon2} Prop 14, and involves integrating along a flow box of the geodesic flow containing $I$.

For the cross-ratio construction, see \cite{Hersonsky}, Prop 4.7. The proof here involves a clever partition of the set of geodesics intersecting $I$, as well as some identities of the cross-ratio function.\qed

\begin{myprop}
For any negatively curved Riemannian surface $X$, the three versions of the Liouville current $L_X$ given here are all equal.
\end{myprop}
{\em Proof.} Let $\alpha\in C$ be an isotopy class of curves on $S$. Recall that the intersection $i(L_X, \alpha)$ is defined to be the $L_X$-measure of the set of all geodesics intersecting a fundamental domain for the action of $\alpha$ on its axis. By the lemma, for any of the three current constructions, this intersection is the length of the fundamental domain, which is the same as the length of the unique geodesic on $X$ in the class of $\alpha$. Thus for any $\alpha$, $i(L_X,\alpha)$ does not depend on which construction of $L_X$ is used. Thus by proposition \ref{otalcurrents}, these currents are equal. \qed 

\renewcommand{\thechapter}{3}

\chapter{Negatively curved cone surfaces}

The addition of cone singularities to a negatively curved surface requires some alterations to the theory described in the previous chapter, but similar results will hold for such surfaces. This chapter does not re-develop the theory in full, but rather outlines the changes from the previous chapter.

A negatively curved cone surface is a surface $Y$ with a negatively curved Riemannian metric which is defined everywhere except at a discrete collection of points called the {\em cone locus} of $Y$. At any cone point, the {\em cone angle} (defined below) is more than $2\pi$ (there are also manifolds with cone points that have cone angles $<2\pi$, but we do not want to consider such surfaces). Note that the cone locus is always finite, since it is discrete and $Y$ is compact. A cone point is also often called a singular point.

For any point $p\in Y$, we define the cone angle at $p$ as follows. For small $\epsilon >0$, let $s_\epsilon$ be the equidistant circle of radius $\epsilon$ at $p$, and $l(s_\epsilon)$ its circumference. Then the cone angle at $p$ is $\lim_{\epsilon\rightarrow 0} l(s_\epsilon)/\epsilon$. If the Riemannian metric is defined on any open neighborhood of $p$, then the cone angle at $p$ is $2\pi$, because Riemannian metrics are infinitesimally Euclidean. At a point where the metric is not defined, the cone angle may be larger. For a constructive way to introduce cone angles into a surface, see \cite{Troyanov1}.

Cone points can be thought of as a way to concentrate some negative curvature into discrete points rather than spreading it out over the surface. Let $Y$ be a negatively curved cone surface with cone locus $P$. Then for each $p_i\in P$, with cone angle $\theta_i$, define $k_i = 2\pi - \theta_i$ to be the {\em concentrated curvature} at $p_i$. There is a version of the Gauss-Bonnet theorem that holds for cone surfaces, which states that
$$\int_{Y\setminus P} K dA + \sum_i k_i = 2\pi\chi(Y),$$
where $K$ is the Gaussian curvature function, and $dA$ is the Riemannian area element.

The following result of Troyanov \cite{Troyanov2} classifies all negatively curved cone surfaces and is an analog of Berger's result in proposition \ref{berger}.

\begin{myprop}
\label{negativeconecurvature}
Let $\Sigma$ be a compact Riemann surface. Choose finitely many points $p_i$ on $\Sigma$ and numbers $\theta_i >0$ so that $\sum (2\pi - \theta_i) > 2\pi\chi(\Sigma)$. Then any smooth negative function on $\Sigma$ is the Gaussian curvature function of a unique cone metric in the conformal class of $\Sigma$, having cone angles $\theta_i$ at $p_i$.
\end{myprop}

\section{Boundary at infinity and conjugacy}

A geodesic on a negatively curved cone structure is defined to be a curve which is a piecewise geodesic in the Riemannian metric away from the cone locus, and forms an angle of at least $\pi$ on both sides at every point along the curve.

Let Neg*($S$) be the moduli space of isotopy classes of negatively curved cone manifolds with markings to $S$, and let $Y$ be such a class of structures. The definition of the boundary at infinity from the page \pageref{boundaryatinfinity} applies without alteration to $Y$, and there is still a one-to-one correspondence between oriented geodesics in $\wt{Y}$ and pairs of distinct points in $\bdyinf Y$.
\label{sectorbehind}
The biggest differences in dealing with cone surfaces as opposed to Riemannian surfaces arise from the following observation. Recall that for a Riemannian surface $X$ and a point $p\in\wt{X}$, the exponential map $T^1_p\wt{X}\rightarrow\bdyinf X$ is a homeomorphism. This is not true for cone surfaces; in fact, there are directions where this map is not even defined.

Let $p\in \wt{Y}$ be a non-singular point, and $V = T^1_p \wt{Y}$ the circle of unit vectors at $p$. A vector $v\in V$ is called a {\em non-singular direction} if the ray from $p$ in the direction of $v$ does not meet any cone points, and a {\em singular direction} otherwise. Let $V_0 \subset V$ be the set of all non-singular directions.  Note that since the finitely-many cone points on $Y$ lift to countably-many cone points on $\wt{Y}$, $V_0$ has full measure in $V$.

The exponential map $V\rightarrow \bdyinf Y$ is only well-defined for the non-singular directions $V_0$. Indeed, let $v$ be a singular direction and $c$ the first cone point reached by traveling from $p$ in the direction of $v$, with cone angle $2\pi + \theta$. Once a geodesic ray reaches $c$, it has a range of $\theta$ possible exit directions, and thus becomes undetermined by $v$. Changing this exit angle will alter the boundary point ultimately reached by the ray. (Of course, if the ray leaves $c$ in a singular direction, it becomes even more undetermined, etc.)

\begin{figure}[htb]
  \centering
  \def\svgwidth{150pt}
  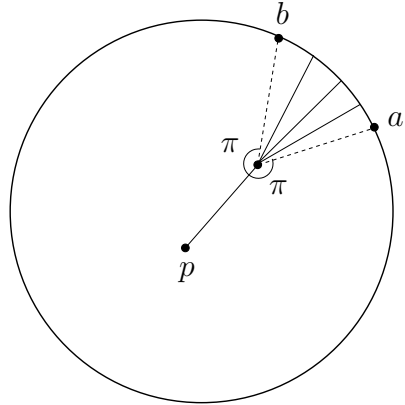
  \caption{The sector behind a cone point.}
\end{figure}

In fact, there is an entire interval $I = [a,b]$ of $\bdyinf Y$ which is ``behind" $c$ from $p$, in the sense that no geodesic ray from $p$ can limit to a point in $I$ without passing through $c$. This interval will be referred to as the {\em sector behind $c$ from $p$}. This argument shows that the complement of the image of $V_0$ in $\bdyinf Y$ contains intervals. Note that since the map $V_0\rightarrow \bdyinf Y$ is increasing, it can be completed to a measurable function on all of $V$ which is either left continuous or right continuous.

It is interesting to note that there is still a one-to-one correspondence between geodesic rays from $p$ and points at infinity (see \cite{Bridson}, Prop 8.2), but such rays are not always determined by their starting directions from $p$.

The universal cover of a negatively curved cone surface is a Gromov-hyperbolic Hadamard space, so proposition \ref{conjugacy} still holds in this setting. In particular, we have:

\begin{myprop}
\label{coneconjugacy}
Given any two classes $Y_1, Y_2 \in \textrm{Neg*}(S)$, there is a unique $\Gamma$-equivariant conjugacy homeomorphism $\bdyinf Y_1 \rightarrow \bdyinf Y_2$. Furthermore, if $X\in \textrm{Neg}(S)$ and $Y\in \textrm{Neg*}(S)$, there is a also a unique $\Gamma$-invariant conjugacy homeomorphism $\bdyinf X\rightarrow\bdyinf Y$.
\end{myprop}

It is important to note here that not only can surfaces be compared within each class of structures, they can also be compared across the two classes.

\section{The Liouville current of a negatively curved cone manifold}

For any $Y \in \textrm{Neg*}(S)$, the space $\mc(Y)$ of (metric) geodesic currents on $Y$ is defined, as in the previous chapter, as the weak* uniform space of all $\pi_1(Y)$-invariant Borel measures on $\mg(\wt{Y}) = (\bdyinf Y \times \bdyinf Y \setminus \Delta)/\mathbf{Z}_2$. Since there is a unique identification $\bdyinf Y \rightarrow \bdyinf S$ via conjugacy, $\mc(Y)$ can again be uniquely identified with the space of (topological) currents $\mc(S)$.

Each $Y\in \textrm{Neg*}(S)$ determines a Liouville current $L_Y\in \mc(S)$ as in the previous chapter, but some care must be taken in the definition to account for the cone points. The three constructions from the previous chapter are discussed below, with these differences noted. In each case there will be differences between the non-singular geodesics (i.e. those which do not meet any cone point) and the singular ones.

As before, all three of these constructions define the same geodesic current, satisfying the property that the measure of the set of all geodesics meeting a given geodesic segment in $\wt{Y}$ is equal to the length of the segment (see section \ref{constructionsequal}).

\subsection{The Liouville current as a transverse measure}
\label{conetransversemeasure}
Let $Y \in \textrm{Neg*}(S)$, with cone locus $P$. Since the Riemannian metric is not defined on $P$, the unit tangent bundle is also not defined above these points. Further, the geodesic flow of $Y$ is only well-defined on the collection of non-singular directions at each non-cone point, since geodesics become undetermined at cone points. Denote the space of non-singular directions on $Y$ as $T_0^1 Y$. Note that $T_0^1 Y$ has full measure in $T^1 Y$ with respect to the Riemannian volume form, since $Y$ has finitely many cone points.

The interior product of the Riemannian volume form with the unit vector field that generates the geodesic flow is then a transverse invariant measure, and this is the Liouville current $L_Y$.

\subsection{The Liouville current from a cross-ratio on $\bdyinf S$}

The definition of the cross-ratio function $C_Y$ on quadruples of distinct points in $\bdyinf S$ goes through unchanged from the previous chapter, as does the construction of the Liouville measure $L_Y$ from $C_Y$. However, $C_Y$ has an important property when $Y$ is a cone surface, which will be important later.

\begin{myprop}
\label{crossratiozero}
Let $\gamma$ be a singular complete geodesic in $\wt{Y}$, with endpoints $\gamma_-, \gamma_+ \in \bdyinf S$. Then there are disjoint, non-trivial intervals $[a,b]$ and $[c,d]$ in $\bdyinf S$ so that $\gamma_-\in [a,b]$, $\gamma_+ \in [c,d]$, and $C_Y(a,b,c,d) = 0$.
\end{myprop}

\begin{figure}[htb]
  \centering
  \def\svgwidth{150pt}
  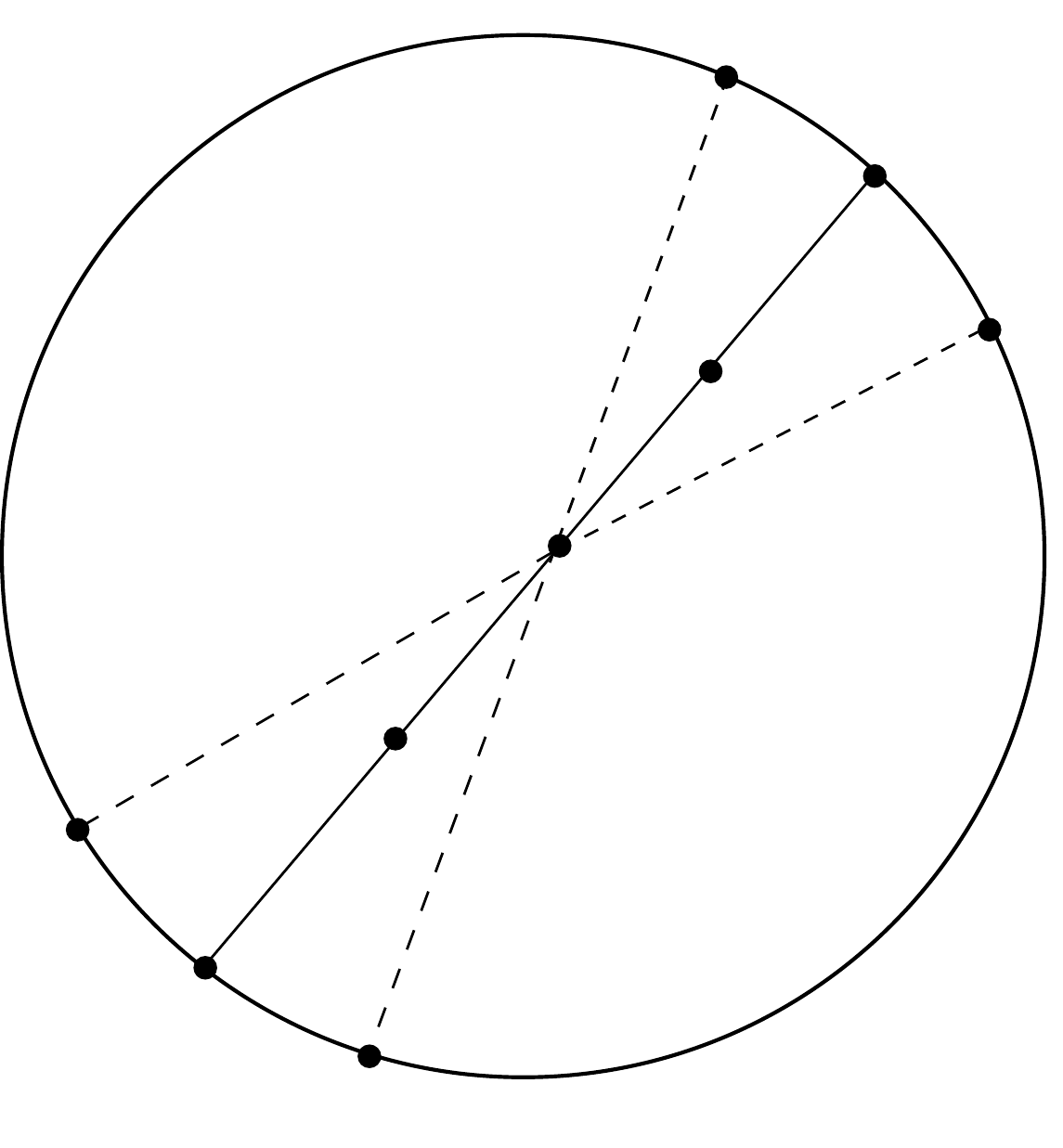
  \caption{Proof of prop \ref{crossratiozero}.}
\end{figure}

{\em Proof.} Let $p$ be a cone point on $\gamma$, and let $x_1$ and $x_2$ be non-cone points on $\gamma$, one on each side of $p$, so that $x_1$ is closer to $\gamma_-$ and $x_2$ is closer to $\gamma_+$, and so that the geodesic segments $\overline{x_1 \ p}$ and $\overline{p \ x_2}$ have no cone points in their interiors. Let $[a,b]$ be the sector behind $p$ from $x_2$ and $[c,d]$ the sector behind $p$ from $x_1$. Then any geodesic with one endpoint in $[a,b]$ and the other in $[c,d]$ must pass through $p$. It is clear from the definition of the cross-ratio on page \pageref{crossratio} that this implies $C_Y(a,b,c,d)=0$. \qed

\subsection{The Liouville current in geodesic-angle coordinates}
\label{conegeodesicangle}
Let $\gamma$ be a non-singular, unit-speed parameterized geodesic on $\wt{Y}$, and let $\mg(\gamma)$ be the set of all geodesics transversely intersecting $\gamma$. As before, we want to create a geodesic-angle coordinate $\xi_{Y,\gamma}:\mg(\gamma)\rightarrow \mathbf{R}\times(0,\pi)$ by sending a geodesic to (the parameter of) its point of intersection with $\gamma$ and the angle of intersection. The problem is that this mapping will not be one-to-one, since two different geodesics can intersect $\gamma$ at the same point and in the same angle if they are both singular.

However, the coordinate becomes a homeomorphism onto its image if we restrict to the subset $\mg_0(\gamma)\subset\mg(\gamma)$ of non-singular geodesics intersecting $\gamma$. Note that the image $\xi_{Y,\gamma}(\mg_0(\gamma))$ has full measure in $\mathbf{R}\times(0,\pi)$ with respect to $d\lambda = 1/2\sin\theta d\theta dt$, since there are only countably many singular directions at each point of $\gamma$.

Let $\mg_0(\wt{Y})\subset \mg(\wt{Y})$ denote the collection of all non-singular geodesics in $\wt{Y}$. Letting $\gamma$ vary, the sets $\mg_0(\gamma)$ form an open cover of $\mg_0(\wt{Y})$, and pulling back the measure $d\lambda$ through each $\xi_{Y,\gamma}$ defines a measure on $\mg_0(\wt{Y})$. Extending this to all of $\mg(\wt{Y})$, by defining the measure to be $0$ outside of $\mg_0(\wt{Y})$, produces the Liouville current $L_Y$.

\subsection{The Liouville current of a cone surface does not have full support}

Recall that the support of a measure is the largest closed set in which every open subset has positive measure. A measure does not have full support if the complement of its support is non-empty, or equivalently if there is at least one open set of measure zero.

\begin{myprop}
\label{conesupport}
For a negatively curved cone surface $Y$ marked to $S$, the Liouville measure $L_Y$ on $\mg(\wt{S})$ does not have full support.
\end{myprop}
{\em Proof.}
From the description of the Liouville measure in geodesic-angle coordinates, it is clear that the support of $L_Y$ is $\overline{\mg_0(\wt{S})}$, the closure of the set of non-singular geodesics. This closure consists of geodesics which are either non-singular themselves, or are a limit of non-singular geodesics.

\begin{figure}[htb]
  \centering
  \def\svgwidth{150pt}
  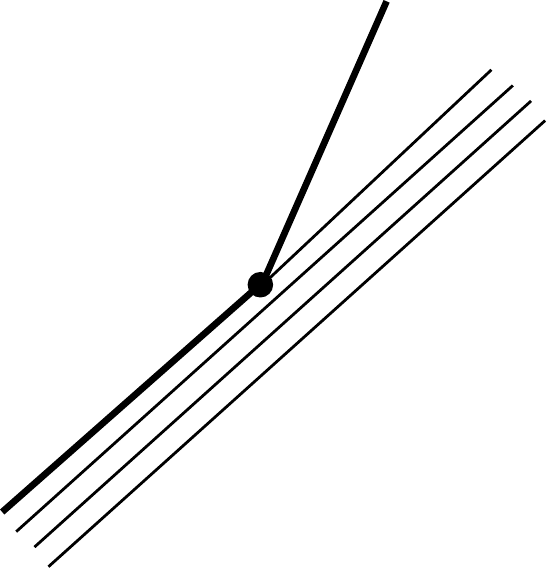
  \caption{A singular geodesic which is not a limit of non-singular geodesics.}
\end{figure}

\label{singularlimit}
There are, however, singular geodesics which are not a limit of non-singular geodesics. Any curve which is a limit of non-singular geodesics must make an angle of exactly $\pi$ on one side at every point along the curve. Therefore any geodesic which passes through a cone point and makes an angle greater than $\pi$ on each side of the singularity cannot be a limit of non-singular geodesics. Thus $\mg(\wt{S}) \setminus \overline{\mg_0(\wt{S})}$ is non-empty, so $L_Y$ does not have full support.

Another argument, based on the cross-ratio, is given by Proposition 2. This says that there exist disjoint, nontrivial intervals $[a,b]$ and $[c,d]$ on $\bdyinf S$ so that $C_Y(a,b,c,d) = 0$. Then the Liouville current satisfies $$L_Y([a,b]\times [c,d]) = |C_Y(a,b,c,d)| = 0$$, so any open set in this rectangle has measure $0$.\qed 

\renewcommand{\thechapter}{4}

\chapter{Flat surfaces}

A flat cone surface is a surface $Z$ with a flat Riemannian metric defined away from a discrete collection of cone points, so that each cone point has a cone angle greater than $2\pi$ (again, there are flat cone surfaces with cone angles less than $2\pi$, but we don't want to consider them). By the Gauss-Bonnet theorem for Riemannian surfaces, any closed flat surface of genus greater than one must necessarily have cone points. Because of this, we can simply say that $Z$ is a {\em flat surface} instead of specifying that $Z$ has cone points.

If cone points are thought of as a way of concentrating negative curvature into discrete points, then a flat surface is one where {\em all} of the negative curvature has been moved into the cone points. The analog of the Gauss-Bonnet theorem for closed flat surfaces says that if $Z$ has cone locus $P = \{p_i\}$ with cone angles $\theta_i$ and concentrated curvatures $k_i = 2\pi - \theta_i$, then
$$\sum_i k_i = 2\pi\chi(Z).$$
Note that since this formula contains no integration against the area element, there is no relationship between the curvature of a flat surface and its area, unlike with negative curvature. This reflects the fact that flat geometries have similarity transformations, whereas curved geometries do not.

The following result of Troyanov \cite{Troyanov1} is an analog of proposition \ref{negativeconecurvature} for flat surfaces.

\begin{myprop}
Let $\Sigma$ be a compact Riemann surface. Choose finitely many points $p_i$ on $\Sigma$ and numbers $\theta_i > 0$ so that $\sum (2\pi - \theta_i) = 2\pi\chi(\Sigma)$. Then for any $A > 0$ there is a unique flat surface of area $A$ in the conformal class of $\Sigma$, having cone points $p_i$ with cone angles $\theta_i$.
\end{myprop}

\label{flatlimit}
A flat surface can be thought of as a limit of negatively curved cone surfaces in the following sense. Let $Z$ be a flat surface, with $p$ a cone point on $Z$ of cone angle $\theta$. For small $\epsilon > 0$, let $Y_\epsilon$ be a negatively curved cone surface with the same cone locus so that:
\begin{itemize}
\item $Y_\epsilon$ has the same cone angles as $Z$, except at $p$, where the angle is $\theta - \epsilon$,
\item away from the cone locus, $Y_\epsilon$ has constant negative Gaussian curvature,
\item $Y_\epsilon$ has the same area $A$ as $Z$,
\item $Y_\epsilon$ is in the same conformal class as $Z$.
\end{itemize}
By proposition \ref{negativeconecurvature}, such a surface must exist, and these requirements determine it uniquely. It is also easy to see that the Gaussian curvature $k_\epsilon$ of $Y_\epsilon$ must be $k_\epsilon = -\epsilon/A$, so as $\epsilon\rightarrow 0$, the curvature goes to 0 and the metrics on $Y_\epsilon$ limit to $Z$. Of course there are many other sequences of negatively curved cone surfaces which limit to $Z$, but this is in some sense the simplest.

\section{Boundary at infinity and conjugacy}

A geodesic on a flat surface is defined to be a curve which is piecewise geodesic away from the cone locus, and which makes an angle of at least $\pi$ on both sides at every point along the curve. The biggest difference between negatively curved cone surfaces and flat surfaces is that in a flat geometry, an equidistant curve to a geodesic is also a geodesic. This is not true in negative (or for that matter, positive) curvature.

Let $\wt{Z}$ be the universal cover of a flat surface, and $\gamma$ a non-singular geodesic in $\wt{Z}$. Then $\gamma$ divides $\wt{Z}$ into two half-spaces, $H_1$ and $H_2$. Let $P_i$ denote the set of cone points contained in $H_i$, and $m_i = \inf_{p\in P_i} d(p,\gamma)$, for $i = 1,2$. If $m_1$ and $m_2$ are both 0, then there are sequences of singularities that limit to $\gamma$ on both sides, and $\gamma$ is in some sense ``trapped" by cone points. If $m_1 > 0$, let $\gamma'$ be a curve in $H_1$ which is equidistant to $\gamma$, with $d(\gamma, \gamma') < m_1$. Let $T$ be the region in $\wt{Z}$ between $\gamma$ and $\gamma'$. Since there are no singularities in $T$, it is isometric to the strip $\mathbf{R}\times [0,d(\gamma,\gamma')]$ in the flat plane, and $\gamma'$ is a geodesic of $\wt{Z}$. This discussion is summarized as a proposition.

\begin{myprop}
\label{stripsexist}
Given a non-singular geodesic $\gamma$ in $\wt{Z}$ which is not limited to by cone points on both sides, $\gamma$ lies in an isometrically embedded strip which is foliated by geodesics equidistant to $\gamma$. The maximal such strip is called the strip of $\gamma$, and $\gamma$ is called a strip geodesic.
\end{myprop}

Now let $Z \in \textrm{Flat*}(S)$ be an isotopy class of flat surfaces. As with negatively curved cone surfaces, the definition of the boundary at infinity $\bdyinf Z$ from page \pageref{boundaryatinfinity} applies to flat surfaces. However, there is no longer a one-to-one correspondence between complete geodesics in $\wt{Z}$ and pairs of distinct points in $\bdyinf Z$. Any geodesic determines a pair of distinct points at infinity, but it is possible for more two different geodesics to determine the same pair of points, for instance if they lie on the same embedded strip. The following proposition says that this is in fact the only way this can happen.

\begin{myprop}
\label{stripisworst}
If two complete geodesics $\gamma_1$ and $\gamma_2$ in $\wt{Z}$ limit to the same pair of points at infinity, then they bound an isometrically embedded strip. In particular, there are no cone points between them.
\end{myprop}
{\em Proof.} Consider first two geodesic rays $\rho_1$ and $\rho_2$ in $\wt{Z}$. There are only three ways in which two rays can limit to the same endpoint at infinity:
\begin{enumerate}
\item Both rays eventually pass through the same cone point and are the same curve afterwards,
\item the rays eventually lie on the same embedded strip, or
\item $\rho_1$ approaches $\rho_2$ by connecting cone points which limit to $\rho_2$.
\end{enumerate}
By the definition of a point at infinity, the complete geodesics $\gamma_1$ and $\gamma_2$ must exhibit one of these behaviors in the positive direction and one in the negative direction. Let i-j denote the case where these two behaviors are i and j.

Case 1-1 violates that $\wt{Z}$ is uniquely geodesic (\cite{Bridson}, prop II.1.4). Cases 1-2 and 1-3 are both impossible by (\cite{Bridson}, prop II.8.2), which says that there is a one-to-one correspondence between $\bdyinf Z$ and the geodesic rays from any point in $\wt{Z}$.

Assume that case 2-3 holds. On the side where the two geodesics lie on a strip, draw a geodesic segment from a point $a$ on $\gamma_1$ to $b$ on $\gamma_2$, perpendicular to both curves. On the side where $\gamma_1$ limits to $\gamma_2$, draw a segment from a cone point $c$ on $\gamma_1$ to $d$ on $\gamma_2$, perpendicular to $\gamma_2$. Then the geodesics and the segments together bound a geodesic polygon in $\wt{Z}$.

\vspace{10mm}
\begin{figure}[htb]
  \centering
  \def\svgwidth{200pt}
  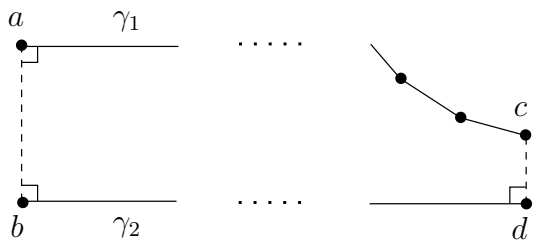
  \caption{An impossible flat geodesic polygon.}
\end{figure}

\noindent Cut this polygon out and double it across its boundary to obtain a flat surface $\Sigma$ homeomorphic to a sphere $S^2$. By the Gauss-Bonnet theorem for flat surfaces, there must be at least a total of $4\pi$ concentrated positive curvature at the cone points of $\Sigma$. The points $a,b,$ and $d$ each contribute $\pi$, and $c$ contributes some amount strictly less than $\pi$. Any cone points in the interior of the polygon can only contribute negative curvature, and the same for any cone points on the boundary of the polygon, since a geodesic makes an angle of at least $\pi$ on both sides when it passes through a cone point. Thus the surface $\Sigma$ can't exist, so this case is impossible.

Case 3-3 is proved similarly to the previous case, and this eliminates all possibilities except for case 2-2. \qed

This shows that there is no canonical way to go from a pair of points on $\bdyinf Z$ to a geodesic in $\wt{Z}$, which will have consequences later. The proposition also implies that, given a curve class $\alpha \in C$, there may not be a unique geodesic on $Z$ in the class of $\alpha$, but any two such will have the same length, since they will be isotopic across an isometrically-embedded flat cylinder. This means that the length spectrum of $Z$ is well-defined as a functional on $C$. Note that this would fail if $Z$ were allowed to have points with cone angle less than $2\pi$.

The discussion on page \pageref{sectorbehind} about singular and non-singular directions applies also to flat surfaces. In particular, the exponential map $T^1_p \wt{Z}\rightarrow \bdyinf Z$ at any non-singular point $p$ is only measurable, and there are sectors on the boundary which are inaccessible from $p$ without passing through cone points.

Since $\wt{Z}$ is a Gromov-hyperbolic Hadamard space, proposition \ref{conjugacy} holds for flat surfaces. Combining this with proposition \ref{coneconjugacy} shows:

\begin{myprop}
\label{bigconjugacy}
Given any two classes $X_1, X_2\in \textrm{NonPos}(S)$, there is a unique $\Gamma$-equivariant conjugacy homeomorphism $\bdyinf X_1\rightarrow \bdyinf X_2$, which can be obtained by lifting the markings to the universal covers and extending to the boundary.
\end{myprop}

\section{The Liouville current of a flat surface}

For any $Z \in \textrm{Flat*}(S)$, the space $\mc(Z)$ of metric geodesic currents on $Z$ is defined as the weak* uniform space of $\pi_1(Z)$-invariant Borel measures on $\mg(\wt{Z}) = (\bdyinf Z \times \bdyinf Z \setminus \Delta)/\mathbf{Z}_2$. Note that, because there may be flat strips in $\wt{Z}$, $\mg(\wt{Z})$ can no longer be precisely identified with the space of geodesics in $\wt{Z}$. The metric currents $\mc(Z)$ can be uniquely identified via conjugacy with the topological currents $\mc(S)$ on $S$.

As before, each $Z \in \textrm{Flat*}(S)$ determines a Liouville current $L_Z \in \mc(S)$. Slightly more care must be taken in defining the Liouville current in the flat case, due to the new phenomenon of flat strips. The construction of $L_Z$ as a transverse measure to the geodesic flow is unchanged from section \ref{conetransversemeasure}. The constructions of $L_Z$ from the cross-ratio and from geodesic-angle coordinates are discussed below.

As in the previous chapters, all three constructions define the same geodesic current, with the property that the measure of the set of all geodesics meeting a given segment in $\wt{Z}$ is equal to the length of the segment (see section \ref{constructionsequal}). Also, the same argument as in proposition \ref{conesupport} shows:

\begin{myprop}
\label{flatconesupport}
For a flat surface $Z$ marked to $S$, the Liouville measure $L_Z$ on $\mg(\wt{S})$ does not have full support.
\end{myprop}

\subsection{The Liouville current from a cross-ratio on $\bdyinf S$}
\label{flatcrossratio}
Recall the definition of the cross-ratio of a surface from section \ref{crossratiodef}, stated here for a flat surface $Z$. Given four distinct points $a,b,c,d\in\bdyinf Z$,
$$C_Z(a,b,c,d) = \frac{1}{2}\lim (d(a_i, c_i) + d(b_i, d_i) - d(a_i, d_i) - d(b_i, c_i)),$$
where the limit is taken as $a_i \rightarrow a$, etc. Here we show that this limit exists, is finite, and does not depend on the choices of sequences tending to the points at infinity.

Given a point $a \in \bdyinf Z$, the {\em Busemann function} $b_a: \wt{Z}\rightarrow\mathbf{R}$ is defined as
$$b_a(x) = \lim_{t\rightarrow\infty} (d(x,\gamma(t)) - t),$$
where $\gamma$ is a geodesic ray limiting to $a$. This limit exists, is finite, and does not depend on choice of $\gamma$ (see \cite{Bridson}, lemma II.8.18). Given two points $x,y \in \wt{Z}$, the {\em horocyclic distance based at $a$} between $x$ and $y$ is $d_a(x, y) = b_a(x) - b_a(y)$. Note that this can be negative, and satisfies the cocycle relation
$$d_a(x,y) + d_a(y,z) = d_a(x,z).$$
Given any $x\in\wt{Z}$, the {\em horocycle based at a} through $x$ is the set of all $y$ so that $d_a(x,y) = 0$. Such a horocycle is a rectifiable curve which limits to $a$ in both directions, and hence bounds a region of $\wt{Z}$ called {\em horoball}. Given a horoball $B$ at $a$, any geodesic limiting to $a$ must eventually lie inside $B$.

Further, given two distinct horocycles $h_1$ and $h_2$ at $a$, and two geodesics $\gamma_1$ and $\gamma_2$ limiting to $a$, the cocycle condition above implies that the segments of $\gamma_1$ and $\gamma_2$ which lie between $h_1$ and $h_2$ have the same length, equal to the absolute value of the horocyclic distance between any point on $h_1$ and any point on $h_2$.

Now let $a,b,c,d \in \bdyinf Z$. Choose disjoint horoballs $H_a, H_b, H_c, H_d$ based at these points. Choose a geodesic connecting $a$ to $c$ (this choice may not be automatic in the flat case by prop \ref{stripsexist}), and let $l(ac)$ be the length of the geodesic segment which lies outside of the horoballs. Similarly define $l(bd), l(ad), l(bc)$. Then following (\cite{Otal2}, lemma 2.1), define
$$\mb(a,b,c,d) = \frac{1}{2} (l(ac) + l(bd) - l(ad) - l(bc)).$$

\begin{figure}[htb]
  \centering
  \def\svgwidth{200pt}
  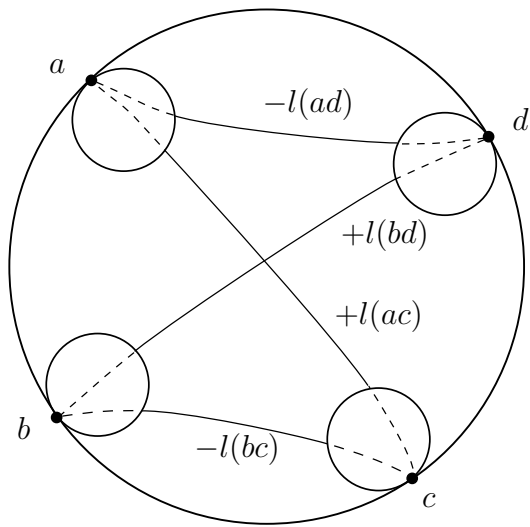
  \caption{The construction of $\mb(a,b,c,d)$.}
\end{figure}

\begin{mylemma}
$\mb(a,b,c,d)$ does not depend on the disjoint horoballs chosen at each point, or on the choice of geodesics connecting $a$ to $c$, etc. Further, $C_Z(a,b,c,d) = \mb(a,b,c,d)$.
\end{mylemma}
{\em Proof.} Given two distinct geodesics $\gamma_1$, $\gamma_2$ connecting $a$ to $c$, they must lie across a flat strip by proposition \ref{stripisworst}. Then the horocycles which bound $H_a$ and $H_c$ cross this flat strip in straight lines perpendicular to the strip. Thus the distance $l(ac)$ does not depend on the choice of geodesic on this strip.

The fact that $\mb$ does not depend on the choice of horoballs follows from the cocycle condition and noticing that each point at infinity has one ``plus" curve and one ``minus" curve limiting to it.

Choose sequences $a_i \rightarrow a$, etc. As $i \rightarrow \infty$, the geodesic segments from $a_i$ to $c_i$ converge to a geodesic from $a$ to $c$. Thus it only remains to show that at $a$ (for instance), the difference of the lengths of the two geodesic segments inside $H_a$ goes to 0 in the limit. This follows because the horocyclic distance based at $a$ between any two points on the same horocycle is 0. \qed

Since $\mb(a,b,c,d)$ is clearly finite and is well-defined by the lemma, so is $C_Z$. As before, the Carath\'eodory construction produces a $\Gamma$-invariant measure on $\mg(\wt{S})$ which uses the cross-ratio to measure rectangles, and this is the Liouville current $L_Z$.

\subsection{The Liouville current in geodesic-angle coordinates}

Let $\gamma$ be a non-singular, unit-speed parameterized geodesic on $\wt{Z}$, and let $\mg(\gamma) \subset \mg(\wt{Z})$ be the collection of all pairs of endpoints of geodesics which intersect $\gamma$ transversely. Note that by prop. \ref{stripsexist}, $\mg(\gamma)$ is not precisely identified with the space of geodesics intersecting $\gamma$.

\begin{figure}[htb]
  \centering
  \def\svgwidth{225pt}
  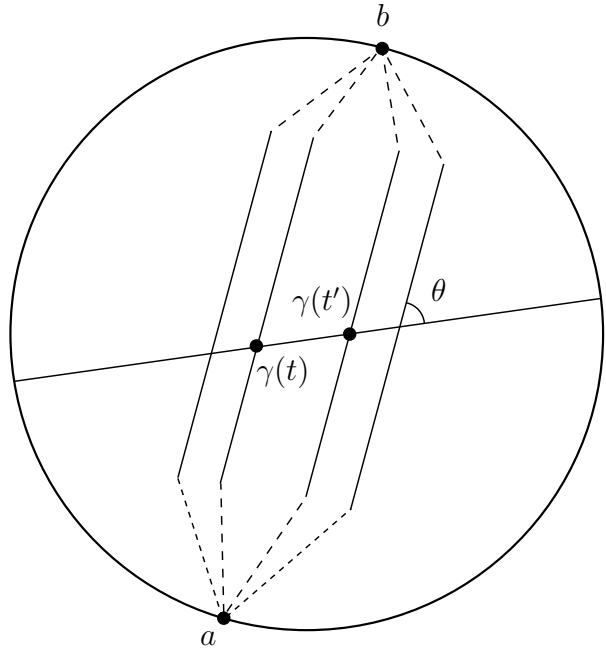
  \caption{Geodesic-angle coordinates are not well-defined across a strip.}
\end{figure}

As in section \ref{conegeodesicangle}, the function $\xi_{Z,\gamma}:\mg(\gamma) \rightarrow \mathbf{R} \times(0,\pi)$, mapping to the intersection point with $\gamma$ and the angle of intersection, is not one-to-one on the set of singular geodesics. Now there is the additional issue that $\xi_{Z,\gamma}$ is also not well-defined if $\gamma$ crosses a strip. Let $a,b \in \bdyinf Z$ be the endpoints of a strip in $\wt{Z}$ which $\gamma$ crosses transversely. Then $\gamma$ crosses this strip at some well-defined angle $\theta$. Choose $t, t'$ so that $\gamma(t), \gamma(t')$ lie on the strip, then $\xi_{Z,\gamma}(a,b)$ could be $(t,\theta)$ or $(t',\theta)$, depending on the geodesic chosen in the strip connecting $a$ and $b$.

The solution to this problem is to think of strip geodesics as being singular for purposes of defining the geodesic-angle coordinates. Let $\mg_0(\gamma)\subset \mg(\gamma)$ be the subset of all endpoints of non-singular, non-strip geodesics intersecting $\gamma$. Then the restriction of $\xi_{Z,\gamma}$ to $\mg_0(\gamma)$ is a homeomorphism onto its image, which has full measure in $\mathbf{R}\times(0,\pi)$ with respect to $d\lambda = 1/2\sin\theta d\theta dt$.

Similarly, let $\mg_0(\wt{Z})\subset\mg(\wt{Z})$ be the collection of all endpoints of non-singular, non-strip geodesics. Then the sets $\mg_0(\gamma)$ form an open cover of $\mg_0(\wt{Z})$, and pulling back $d\lambda$ through each $\xi_{Y,\gamma}$ defines a measure on $\mg_0(\wt{Z})$. Extend this to $\mg(\wt{Z})$ by defining the measure to be 0 outside of $\mg_0(\wt{Z})$ to obtain the Liouville current $L_Z$. 
\renewcommand{\thechapter}{5}

\chapter{Spectral rigidity}

Recall that $S$ is a topological surface of genus $\geq 2$, and $C$ is the set of isotopy classes of closed curves on $S$.

Given any non-positively curved surface $X$ with a marking to $S$, the length spectrum of $X$ is a function $l_X: C\rightarrow \mathbf{R}$ which takes a curve class $\alpha$ on $S$ to the length of the geodesic on $X$ in the class of $\alpha$. If $X$ is strictly negatively curved this geodesic is unique; if $X$ is flat then there may be more than one geodesic in the class of $\alpha$, but any two will have the same length since they will be isotopic across an embedded cylinder. Thus the length spectrum is always well defined in $\mathbf{R}^C$.

It is clear that if two isometric surfaces are marked to $S$ via isotopic markings, they define the same length spectrum. Thus if we define NonPos($S$) to be the union of Neg($S$), Neg*($S$), and Flat*($S$), the length spectrum defines a function $l: \textrm{NonPos}(S)\rightarrow \mathbf{R}^C$. There is also a function $L: \textrm{NonPos}(S)\rightarrow \mc(S)$ that takes any surface to its Liouville measure, and on page \pageref{otalcurrents}, we defined a function $I: \mc(S)\rightarrow \mathbf{R}^C$ that maps $\mu \mapsto i(\mu,-)$.

\vspace{10mm}

\begin{center}
\begin{tikzpicture}[description/.style={fill=white,inner sep=2pt}]
\matrix (m) [matrix of math nodes, row sep=3em,
column sep=2.5em, text height=1.5ex, text depth=0.25ex]
{ \textrm{NonPos}(S) & & \rc \\
& \mc(S) & \\ };
\path[->,font=\scriptsize]
(m-1-1) edge node[auto] {$ l $} (m-1-3)
edge node[auto] {$ L $} (m-2-2)
(m-2-2) edge node[auto] {$ I $} (m-1-3);
\end{tikzpicture}
\end{center}

\begin{myprop}
\label{crofton}
The above diagram is commutative. That is, $l = I \circ L$ as functions $\textrm{NonPos}(S)\rightarrow \mathbf{R}^C$.
\end{myprop}
{\em Proof.} Let $X\in\textrm{NonPos}(S)$, and $\alpha\in C$. Then $I(L_X)(\alpha) = i(L_X, \alpha)$, so it is equivalent to prove that $i(L_X, \alpha) = l_X(\alpha)$ for all $\alpha$. That is, that the intersection of $L_X$ with any curve class is equal to the length of the $X$-geodesic in that class.

Think of $\alpha$ as a conjugacy class of $\Gamma$, and let $\gamma\in\alpha$ be a specific representative. Let $J$ be a fundamental domain for the action of $\gamma$ on its axis in $\wt{X}$. Then $l_X(\alpha)$ is equal to the length of $J$. Let $\mg(J)$ be the set of all geodesics intersecting $J$.

Since $J$ is a geodesic segment in $\wt{X}$, its length is equal to $L_X(\mg(J))$, by lemma \ref{lengthismeasure}, which also holds for cone surfaces. But this is $i(L_X,\alpha)$, by the definition of the intersection pairing. Thus $l_X(\alpha) = L_X(\mg(J)) = i(L_X,\alpha)$, and the result is shown. \qed

\section{Length spectrum separates cone surfaces from Riemannian surfaces}

In this section it is shown that no negatively curved Riemannian surface has the same length spectrum as either a negatively curved cone surface or a flat surface. That is, if we let NonPos*($S$) denote the union of Neg*($S$) and Flat*($S$), then the images $l(\textrm{Neg}(S))$ and $l(\textrm{NonPos*}(S))$ do not overlap in $\mathbf{R}^C$. This only requires the equivalence of the first two items in the following proposition, but the others are recorded here for future use.

\begin{myprop}
\label{spectrumequiv}
Let $X_1, X_2 \in \textrm{NonPos}(S)$, and let $\phi : \bdyinf X_1 \rightarrow \bdyinf X_2$ be the conjugacy map between them. Then the following are equivalent:
\begin{enumerate}[(1)]
\item $X_1$ and $X_2$ have the same length spectrum in $\mathbf{R}^C$,
\item $X_1$ and $X_2$ have the same Liouville current in $\mc(S)$,
\item $\phi$ takes the cross-ratio of $X_1$ to the cross-ratio of $X_2$, so $C_{X_1}(a,b,c,d) = \\ C_{X_2}(\phi(a),\phi(b),\phi(c),\phi(d))$,
\item Let $\gamma$ be a geodesic on $\wt{X_1}$ with endpoints $a,b$ at infinity, and $\gamma'$ the associated geodesic in $\wt{X_2}$ with endpoints $\phi(a),\phi(b)$. Then the map $\phi\times\phi$ is a measure-isomorphism of $d\lambda$ in the geodesic-angle coordinates on $\mg(\gamma)$ and $\mg(\gamma')$.
\end{enumerate}
\end{myprop}
{\em Proof.} $(2)\iff(3)$ is obvious, by the construction of the Liouville currents from the cross-ratios. Similarly, $(2)\iff(4)$ follows from the construction of the Liouville currents from geodesic-angle coordinates. It remains to show $(1)\iff(2)$.

Assume $l_{X_1} = l_{X_2}$. Then $I(L_{X_1}) = I(L_{X_2})$ by Prop \ref{crofton}, and $I$ is injective by proposition 1.2. Thus $L_{X_1} = L_{X_2}$.

Conversely, assume $L_{X_1} = L_{X_2}$. Then $I(L_{X_1}) = I(L_{X_2})$, and $I\circ L = l$ by Prop \ref{crofton}, so $l_{X_1} = l_{X_2}$. This completes the proof. \qed

\begin{mythm}
Let $X$ be a negatively curved Riemannian surface and $Y$ either a negatively curved cone surface or a flat surface, each with a marking to $S$. Then $X$ and $Y$ have different length spectra.
\end{mythm}
{\em Proof.} Let $U \subset \mg(\wt{S})$ be a small open set, contained in the domain of some geodesic-angle coordinate $\xi_{X,\gamma}: \mg(\gamma)\rightarrow \mathbf{R}\times (0,\pi)$. Then $L_X(U) = \int_{\xi_{X,\gamma}(U)}1/2\sin\theta d\theta dt$ is positive, since $\xi_{X,\gamma}$ is a homeomorphism. This shows that $L_X$ has full support on $\mg(\wt{S})$.

On the other hand, propositions \ref{conesupport} and \ref{flatconesupport} state that $L_Y$ does not have full support, so $L_X \neq L_Y$ as currents in $\mc(S)$. By prop \ref{spectrumequiv}, this implies that $X$ and $Y$ have different length spectra. \qed

\section{Rigidity for negatively curved Riemannian surfaces}
\label{otalsproof}
In this section we present an outline of Otal's argument in \cite{Otal} that Neg($S$) is spectrally rigid, i.e. no two distinct classes in Neg($S$) determine the same length spectrum.

Let $X_1, X_2$ be negatively curved Riemannian surfaces, with markings $f_i:X_i\rightarrow S$. It will be shown that if $X_1$ and $X_2$ have the same length spectrum in $\mathbf{R}^C$, then there is an isometry $h:X_1\rightarrow X_2$ which is isotopic to $f_2^{-1}\circ f_1$. This implies that $X_1$ and $X_2$ belong to the same class of surfaces in Neg($S$).

Let $\phi:\bdyinf X_1 \rightarrow\bdyinf X_2$ be the conjugacy map between the two surfaces, and $\Phi =\phi\times\phi:\mg(\wt{X_1})\rightarrow\mg(\wt{X_2})$ the associated correspondence between their spaces of geodesics. Given two intersecting geodesics $\alpha$, $\beta$ in $X_1$, the corresponding geodesics $\Phi(\alpha), \Phi(\beta)$ must intersect in $\wt{X_2}$, since two geodesics intersect if and only if their endpoints are interlaced on the boundary.

However, given three geodesics $\alpha, \beta, \gamma$ which all pass through a common point in $\wt{X_1}$, it is not necessary that $\Phi(\alpha), \Phi(\beta), \Phi(\gamma)$ all share a common point in $\wt{X_2}$. In general they will form a geodesic triangle, denoted $T(\alpha,\beta,\gamma)$. Otal's method was to show that if the two marked surfaces have the same length spectrum, then the sum of the interior angles of $T(\alpha,\beta,\gamma)$ is $\pi$, which implies by the Gauss-Bonnet formula (see \cite{Lee}, p. 164) that the triangle degenerates to a single point, since the curvature is strictly negative. This implies that the collection of all geodesics passing through any point $p$ in $\wt{X_1}$ maps under $\Phi$ to the collection of all geodesics passing through some point $p'$ in $\wt{X_2}$. Then the isometry $h$ is defined by sending $p$ to $p'$.

To begin to show this, define a continuous function $\theta': T^1 X_1\times [0,\pi]\rightarrow [0,\pi]$ as follows. Given a unit vector $v$ based at $p\in X_1$ and $\theta\in [0,\pi]$, let $\gamma_v$ be the geodesic through a lift $\wt{p}\in\wt{X_1}$ of $p$ in the direction of $v$, and $\gamma_{\theta\cdot v}$ the geodesic through $\wt{p}$ in the direction of $\theta\cdot v$. Then $\theta'(v,\theta)$ is the angle in $\wt{X_2}$ between $\Phi(\gamma_v)$ and $\Phi(\gamma_{\theta\cdot v})$. Define $\Theta':[0,\pi]\rightarrow [0,\pi]$ so that $\Theta'(\theta)$ is the average of $\theta'(v,\theta)$ over $v\in T^1 X_1$. In other words, $\Theta'(\theta)$ is the average angle between pairs of geodesics in $X_2$ which meet at angle $\theta$ in $X_1$.

\begin{mylemma}
$\Theta'$ is an increasing homeomorphism of $[0,\pi]$ satisfying:
\begin{itemize}
\item $\Theta'$ is symmetric in $\pi - \theta$, so $\Theta'(\pi-\theta) = \pi - \Theta'(\theta)$. Note that this is equivalent to the graph being rotationally symmetric about the midpoint $(\pi/2,\pi/2)$.
\item $\Theta'$ is super-additive, so $\Theta'(\theta_1 + \theta_2) \geq \Theta'(\theta_1) + \Theta'(\theta_2)$ when $\theta_1 + \theta_2 \leq \pi$.
\end{itemize}
\end{mylemma}
\label{superadditive}
{\em Proof.} Only the super-additivity will be outlined here, since it will be needed later. Choose three geodesics $\gamma_v, \gamma_{\theta_1\cdot v}, \gamma_{(\theta_1+\theta_2)\cdot v}$ which pass through a common point in $\wt{X_1}$. Then the interior angles of the geodesic triangle $T(\gamma_v, \gamma_{\theta_1\cdot v}, \gamma_{(\theta_1+\theta_2)\cdot v})$ in $\wt{X_2}$ are $\theta'(v,\theta_1)$, $\theta'(\theta_1\cdot v,\theta_2)$, and $\pi - \theta'(v, \theta_1+\theta_2)$. Since $\wt{X_2}$ is negatively curved, the sum of these angles must not exceed $\pi$, so:
$$\theta'(v,\theta_1) + \theta'(\theta_1\cdot v,\theta_2) \leq \theta'(v,\theta_1 +\theta_2).$$
Integrating this expression first over the fiber of $T^1 X_1$ over $p$ (using the fact that the Lebesgue measure on this fiber is invariant under rotation), and then integrating again over the surface, gives the result. \qed

\begin{myprop}
\label{thetaprimeisidentity}
If $X_1$ and $X_2$ have the same length spectrum, $\Theta'$ is the identity map.
\end{myprop}
{\em Proof.} Let $F:[0,\pi]\rightarrow\mathbf{R}$ be continuous and convex. Then by Jensen's inequality \cite{Rudin}, we have
$$F(\Theta'(\theta)) \leq \frac{1}{V(T^1 X_1)}\int_{T^1 X_1} F(\theta'(v,\theta))dv,$$
where $dv$ is the volume form on $T^1 X_1$, invariant under the geodesic flow. Integrate both sides with respect to the measure $\sin\theta d\theta$ on $[0,\pi]$ and switch the order of integration on the right:
$$\int_0^\pi F(\Theta'(\theta))\sin\theta d\theta \leq \frac{1}{V(T^1 X_1)}\int_{T^1 X_1}\int_0^\pi F(\theta'(v,\theta))\sin\theta d\theta dv.$$
Define $F'(v) = \int_0^\pi F(\theta'(v,\theta))\sin\theta d\theta$. Then the right side of the above inequality is the average of $F'$ over $T^1 X_1$. A lemma is now needed.

\begin{mylemma}
The average of $F'$ on $T^1 X_1$ is equal to $\int_0^\pi F(\theta)\sin\theta d\theta$.
\end{mylemma}
{\em Proof.} Since $dv$ is a measure on $T_1 X^1$ invariant under the geodesic flow, it is a limit of measures supported on closed geodesics (see page \pageref{curvesdenseincurrents} and \cite{BonahonThesis}). Let $\gamma$ be a closed orbit of the geodesic flow, then the average of $F'$ on $\gamma$ is
\begin{eqnarray*}
 \frac{1}{l(\gamma)}\int_\gamma F'(\gamma(t))dt &=& \frac{1}{l(\gamma)}\int_{\gamma\times (0,\pi)} F(\theta'(\gamma(t),\theta))\sin\theta d\theta dt \\ \\
   &=& \frac{2}{l(\gamma)}\int_{\mg(\gamma)} F(\theta'(\gamma(t),\theta))d\lambda,
\end{eqnarray*}
where $\mg(\gamma)$ is the set of all geodesics intersecting $\gamma$ and $d\lambda = 1/2\sin\theta d\theta dt$, as in previous chapters. Let $\gamma'$ be the associated geodesic in $X_2$ to $\gamma$. Since $X_1$ and $X_2$ have the same length spectrum, proposition \ref{spectrumequiv} says that $\Phi: \mg(\gamma)\rightarrow\mg(\gamma')$ takes $d\lambda$ to $d\lambda' = 1/2\sin\theta' d\theta' dt'$. Then change variables via $\Phi$ to get
\begin{eqnarray*}
\frac{2}{l(\gamma)}\int_{\mg(\gamma)} F(\theta'(\gamma(t),\theta))d\lambda &=& \frac{2}{l(\gamma)} \int_{\mg(\gamma')} F(\theta')d\lambda' \\ \\
 &=& \frac{l(\gamma')}{l(\gamma)}\int_0^\pi F(\theta')\sin\theta' d\theta' \\ \\
 &=& \int_0^\pi F(\theta)\sin\theta d\theta,
\end{eqnarray*}
since $l(\gamma) = l(\gamma')$ by assumption. This shows that the average of $F'$ along any closed geodesic is $\int_0^\pi F(\theta)\sin\theta d\theta$, so this must also be the average over $T^1 X_1$, since $dv$ is a limit of measures supported on closed geodesics.

\noindent{\em Back to prop \ref{thetaprimeisidentity}.} Apply the above lemma to the last inequality to obtain \label{integralinequality}
$$\int_0^\pi F(\Theta'(\theta))\sin\theta d\theta \leq \int_0^\pi F(\theta)\sin\theta d\theta.$$
Otal then proves a lemma that given an increasing, super-additive homeomorphism $\Psi$ of $[0,\pi]$ to itself which is symmetric in $\pi - \theta$ and satisfies
$$\int_0^\pi F(\Psi(\theta))\sin\theta d\theta \leq \int_0^\pi F(\theta)\sin\theta d\theta$$
\noindent for any convex $F$, $\Psi$ must be the identity (see lemma \ref{otallemma} below for a more general proof of this fact where $\Psi$ need only be measurable). Thus $\Theta'$ is the identity. \qed

Since $\Theta'$ is the identity map, it is in particular not just super-additive, but strictly additive, so $\Theta'(\theta_1 + \theta_2) = \Theta'(\theta_1) + \Theta'(\theta_2)$. Looking back at the argument on page \pageref{superadditive} that $\Theta'$ is super-additive, this implies that given three geodesics $\alpha,\beta,\gamma$ through a common point $p \in \wt{X_1}$, the triangle $T = T(\alpha,\beta,\gamma)$ has interior angles summing to $\pi$, so $T$ must degenerate to a single point $p'$. As suggested above, define a map $h: \wt{X_1}\rightarrow\wt{X_2}$ so that $p \mapsto p'$.

Let $p,q \in \wt{X_1}$, and $I$ the geodesic segment connecting them. Let $I'$ be the geodesic segment connecting $h(p), h(q)$ in $\wt{X_2}$. Then by lemma \ref{lengthismeasure} and prop \ref{spectrumequiv},
$$d(p,q) = L_{X_1}(\mg(I)) = L_{X_2}(\mg(I')) = d(h(p),h(q)),$$
so $h$ is an isometry. Since for any $\gamma \in \Gamma$, $h$ takes the axis of $\gamma$ in $\wt{X_1}$ to the axis of $\gamma$ in $\wt{X_2}$, we get that $h$ is in the same mapping class as $f_2^{-1}\circ f_1$, so $X_1$ and $X_2$ with their markings to $S$ are in the same class in Neg($S$). This completes the proof.

\section{Rigidity for cone surfaces}

Hersonsky and Paulin in \cite{Hersonsky} adapted Otal's proof from the previous section to show that Neg*($S$) is spectrally rigid. The important observation is that the methods used do not require that the functions defined in the proof ($\theta'$, $\Theta'$, $F'$, etc) are continuous, only that they are measurable.

Let $f_i:Y_i\rightarrow S$, $i = 1,2$, be marked negatively curved cone surfaces having the same length spectrum in $\mathbf{R}^C$. Two geodesics $\alpha,\beta\in\mg(\wt{Y_1})$ may intersect in a geodesic segment, but if they are chosen to be non-singular, this will not happen.

Let $T_0 \subset T^1 X_1 \times [0,\pi]$ be all $(v,\theta)$ such that both $v$ and $\theta\cdot v$ are non-singular directions. Note that this is a subset of full measure for the product of the volume on $T^1 X_1$ and the Lebesgue measure on $[0,\pi]$. As before, for any $(v,\theta) \in T_0$, let $\gamma_v$ and $\gamma_{\theta\cdot v}$ be geodesics in $\wt{Y_1}$ determined by lifts of $v$ and $\theta \cdot v$. The goal is to define $\theta'(v,\theta)$ to be the angle at which $\Phi(\gamma_v)$ and $\Phi(\gamma_{\theta\cdot v})$ intersect, but first we must know that these two geodesics are also non-singular. But this is clear because $\Phi$ takes the support of $L_{Y_1}$ to the support of $L_{Y_2}$ (see propositions \ref{spectrumequiv} and \ref{conesupport}).

Then $\Theta'(\theta)$, defined as before as the average of $\theta'(v,\theta)$, is not continuous, but it is measurable and increasing. The rest of the proof goes through without alteration, and shows that the collection of all non-singular geodesics through any non-cone point $p\in\wt{Y_1}$ is mapped via $\Phi$ to the collection of all non-singular geodesics through some non-cone point $p'\in\wt{Y_2}$. This defines the isometry $h$ away from the cone points, so there is a unique extension to an isometry $h:\wt{Y_1}\rightarrow\wt{Y_2}$.

The author has adapted this argument to show the following result.

\begin{mythm}
\label{conetheorem}
The images of the length spectrum mappings $l: \textrm{Neg*}(S)\rightarrow\mathbf{R}$ and $l: \textrm{Flat*}(S)\rightarrow\mathbf{R}$ do not overlap. That is, no negatively curved cone surface has the same length spectrum as a flat surface.
\end{mythm}
{\em Proof.} Let $Y$ be a marked negatively curved cone surface, and $Z$ a marked flat surface. Assume by absurd that their length spectra are the same.

As above, define $T_0 \subset T^1 Z\times [0,\pi]$ to be all $(v,\theta)$ so that $v$ and $\theta\cdot v$ are both non-singular directions, and $\theta'(v,\theta)$ to be the angle in $\wt{Y}$ between $\Phi(\gamma_v)$ and $\Phi(\gamma_{\theta\cdot v})$. Let $\Theta'(\theta)$ be the average of $\theta'(v,\theta)$ over all $v\in T^1 Z$ so that $(v,\theta)\in T_0$.

Then $\Theta':[0,\pi]\rightarrow [0,\pi]$ is measurable and increasing, and Otal's argument that $\Theta'$ is symmetric in $\pi - \theta$ and super-additive applies without alteration, as does the proof that for any convex $F$,
$$\int_0^\pi F(\Theta'(\theta))\sin\theta d\theta \leq \int_0^\pi F(\theta)\sin\theta d\theta$$
\noindent (see page \pageref{integralinequality}). Then the following lemma is a measurable extension of lemma 8 in \cite{Otal}:

\begin{mylemma}
\label{otallemma}
Let $\Psi$ be an increasing measurable function from $[0,\pi]$ to itself so that
\begin{itemize}
\item $\Psi$ is super-additive and symmetric in $\pi - \theta$,
\item For any convex function $F$ on $[0,\pi]$, ${\displaystyle \int_0^\pi F(\Psi(\theta))\sin\theta d\theta \leq \int_0^\pi F(\theta)\sin\theta d\theta}$
\end{itemize}
\noindent Then $\Psi$ is the identity.
\end{mylemma}
{\em Proof.} There is no interval $(0,a)$ on which $\Psi(x) < x$. If there were, the convex function $F_a(x) = \max(a-x,0)$ would contradict the second condition above. There is then a sequence $\{x_i\}$ limiting to 0 so that $\Psi(x_i) \geq x_i$.

Suppose $\Psi$ is not the identity. Choose $b$ so that $\Psi(b)\neq b$, and since $\Psi$ is symmetric in $\pi - \theta$ (i.e. turn the graph over symmetry), assume that $\Psi(b) < b$. Let $a = \sup\{x<b : \Psi(x)>x\}$. This $\sup$ is not taken over an empty set by the choice of $\{x_i\}$.

Let $c \in (\Psi(b), b)$. Then $\Psi(c) < \Psi(b)$ since $\Psi$ is increasing, and $\Psi(b) < c$. Thus $\Psi(c) < c$ for every $c \in (\Psi(b),b)$. This implies that $a \leq \Psi(b) < b$ by definition of $a$.

Assume that $\Psi(a) > a$. Then $a < \Psi(a) < \Psi(b) < b$, since $a < b$ and $\Psi$ is increasing. Thus $\Psi(a)$ is between $a$ and $b$, so $\Psi(\Psi(a)) < \Psi(a)$ by definition of $a$. But this contradicts that $\Psi$ is increasing, since by assumption $\Psi(a) > a$.

Now assume $\Psi(a) < a$, and let $c \in (\Psi(a), a)$. Then $\Psi(c) < \Psi(a)$ by increasing, and $\Psi(a) < c$, so $\Psi(c) < c$ for all $c\in (\Psi(a),a)$. This contradicts the definition of $a$.

So $\Psi(a) = a$. Now choose $i$ large enough so that $a + x_i \in [a,b]$. Then
$$a + \Psi(x_i) = \Psi(a) + \Psi(x_i) \leq \Psi(a+x_i) < a + x_i,$$
\noindent where the middle inequality is the super-additivity of $\Psi$, and the last follows from the definition of $a$. Then $\Psi(x_i) < x_i$, which contradicts the definition of $x_i$. Therefore $\Psi$ is the identity. \qed

\noindent {\em Back to Theorem \ref{conetheorem}.} By the lemma, $\Theta'$ is the identity on $[0,\pi]$. The lemma is remarkable because it takes a function which is initially assumed only to be measurable and proves that it is the identity. Then the same argument as in the previous section shows that given three co-incident non-singular geodesics in $\wt{Z}$, the geodesic triangle obtained in $\wt{Y}$ by conjugating them has interior angles that add up to $\pi$. Since the curvature of $\wt{Y}$ is strictly negative, this triangle must degenerate to a single point.

As before, this allows us to create an isometry from $\wt{Z}$ to $\wt{Y}$, by conjugating the sheaves of geodesics which pass through each non-singular point. But this is clearly a contradiction, since $Z$ is flat and $Y$ is negatively curved and hence they are not isometric. Thus the marked surfaces do not have the same length spectrum. \qed

Note that this proof would have broken down only at the very end if the roles of the two surfaces were interchanged. If three geodesics through a point in $\wt{Y}$ are conjugated to $\wt{Z}$ and the interior angles of the resulting triangle add up to $\pi$, this only implies that the triangle does not bound any singularities, since $\wt{Z}$ is flat. Keeping the negatively curved surface as the surface which is conjugated {\em to} allows us to use the favorable Gauss-Bonnet formula and complete the proof.

\section{Rigidity for flat surfaces}
\label{flatrigidity}

Let $Z_1$ and $Z_2$ be flat surfaces with markings to $S$, and let $\phi$ be the conjugacy map between their boundaries. One would like to prove that if these surfaces define the same length spectrum in $\rc$, then there is an isometry $h: Z_1 \rightarrow Z_2$ isotopic to the composition of the markings. As noted at the end of the previous section, Otal's proof does not adapt directly to this case, because there are non-trivial flat triangles whose angles sum to $\pi$. We outline here an sketch of an incomplete possible proof.

Assume that the two flat surfaces have the same length spectrum. Let $Y_t$, for $t \in (0,1]$, be a deformation of negatively curved cone surfaces which limits to $Z_2$ as $t\rightarrow 0$. This deformation can be chosen (as described on page \pageref{flatlimit}) so that the cone locus does not change, only one cone angle is altered, $Y_t$ has constant curvature $-t$ away from the cone points, and the area and conformal class do not change. Let $\phi_t$ denote the conjugacy map $\bdyinf Z_1 \rightarrow \bdyinf Y_t$.

\begin{myprop}
\label{spectraconverge}
As $t \rightarrow 0$, the length spectra $l_{Y_t}$ converge in $\rc$ to $l_{Z_2}$, and therefore the Liouville currents converge also in $\mc(S)$.
\end{myprop}
{\em Proof.} Note that topological convergence in $\rc$ is equivalent to pointwise convergence of functionals, so it is equivalent to show that for each $\alpha \in C$, $l_{Y_t}(\alpha) \rightarrow l_{Z_2}(\alpha)$.

Let $\gamma$ be a $Z_2$-geodesic in the class of $\alpha$, and let $\gamma_t$ be $Y_t$-geodesics in the class of $\alpha$. Clearly $l_{Z_2}(\gamma) < l_{Z_2}(\gamma_t)$, and similarly $l_{Y_t}(\gamma_t) < l_{Y_t}(\gamma)$. By the choice of the deformation, $l_{Z_2}(\gamma_t) < l_{Y_t}(\gamma_t)$, so altogether,
$$l_{Z_2}(\gamma) < l_{Z_2}(\gamma_t) < l_{Y_t}(\gamma_t) < l_{Y_t}(\gamma).$$
As $t\rightarrow 0$, $l_{Y_t}(\gamma) \rightarrow l_{Z_2}(\gamma)$, and $l_{Y_t}(\gamma_t)$ is squeezed between these two, so the length spectrum converges. The Liouville currents then also converge, by the injectivity of the map $I: \mc(S)\rightarrow\rc$ and the completeness of $\mc(S)$. \qed

\vspace{10mm}
\begin{figure}[htb]
  \centering
  \def\svgwidth{200pt}
  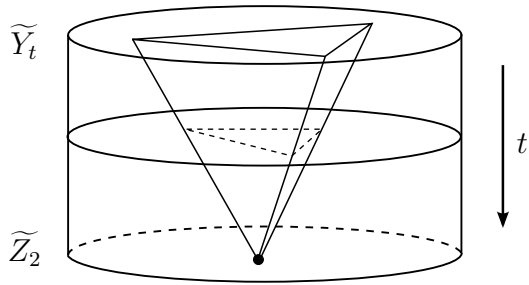
  \caption{A portion of the universal covers of the surfaces in the deformation.}
\end{figure}

The gist of the potential proof is as follows. Since $Z_1$ is not isometric to any $Y_t$, there are three geodesics in $\wt{Z_1}$ which pass through a common point, but which create a non-trivial triangle when conjugated to $Y_t$. Using Otal's average angle function, show that the sum of the angles of such a triangle goes to $\pi$ as $t \rightarrow 0$, so the triangles shrink to a single point in $\wt{Z_2}$. Then it must be shown that the limit of $\wt{Y_t}$-geodesics, with points at infinity fixed, is a $Z_2$-geodesic. This would show that the three co-incident $\wt{Z_1}$-geodesics conjugate to three co-incident $\wt{Z_2}$ geodesics, and thus an isometry can be constructed as in Otal's proof.

Let $T_t \subset T^1(Z_1)\times [0,\pi]$ be the set of all $(v,\theta)$ so that $\gamma_v, \gamma_{\theta\cdot v}, \phi_t(\gamma_v), \phi_t(\gamma_{\theta\cdot v})$ are all non-singular. Note that this is a subset of full measure. Define $\theta_t': T_t\rightarrow [0,\pi]$ so that $\theta_t'(t,\theta)$ is the angle between $\phi_t(\gamma_v)$ and $\phi_t(\gamma_{\theta\cdot v})$ in $\wt{Y_t}$. Then set
$$\Theta_t'(\theta) = \frac{1}{V(T^1 Z_1)} \int_{T^1 Z_1} \theta_t'(v,\theta) dv.$$
As in the previous sections, $\Theta_t': [0,\pi] \rightarrow [0,\pi]$ is increasing, super-additive, and symmetric in $\pi - \theta$.

The following proposition is a generalization of prop \ref{thetaprimeisidentity}, which says that when two surfaces have the same spectra, the average angle function between them is the identity. This new proposition shows that as $t\rightarrow 0$, the average angle functions $\Theta_t$ on $[0,\pi]$ go to the identity.

\begin{myprop}
For any small $t$, there is $\ep(t) > 0$ so that $\sup_{\theta\in [0,\pi]}|\Theta_t'(\theta) - \theta| < \ep$, and this can be chosen so that $\ep \rightarrow 0$ as $t\rightarrow 0$.
\end{myprop}

For any $a \in [0,\pi]$, define the convex function $F_a(\theta) = \max (a-\theta,0)$. Jensen's inequality implies that for any $a$,
$$F_a(\Theta_t'(\theta)) \leq \frac{1}{V(T^1 Z_1)} \int_{T^1 Z_1} F_a(\theta_t'(v,\theta))dv.$$
Then integrate this inequality over $[0,\pi]$ with respect to $\sin\theta d\theta$ and exchange integrals on the right:
$$\int_0^\pi F_a(\Theta_t'(\theta))\sin\theta d\theta \leq \frac{1}{V(T^1 Z_1)} \int_{T^1 Z_1} \left( \int_0^\pi F_a(\theta_t'(v,\theta))\sin\theta d\theta\right)dv.$$

Now define $F_a' = \int_0^\pi F_a(\theta_t'(v,\theta))\sin\theta d\theta$. Then the right side of the above inequality is the average of $F_a'$ over $T^1 Z_1$, with respect to the invariant measure $dv$. As in section \ref{otalsproof}, $dv$ is the limit of measures supported along single orbits of the geodesic flow, so choose a geodesic $\gamma$ on $Z_1$ and average $F_a'$ over $\gamma$:
\begin{eqnarray*}
\frac{1}{l_{Z_1}(\gamma)}\int_\gamma \int_0^\pi F_a(\theta_t'(v,\theta))\sin\theta d\theta dt &=& \frac{2}{l_{Z_1}(\gamma)} \int_{\mg(\gamma)} F_a(\theta_t'(v,\theta)) d\lambda \\ \\
 &=& \frac{2}{l_{Z_1}(\gamma)} \int_{\mg(\gamma')} F_a(\theta')(\Phi_t^* d\lambda)
\end{eqnarray*}
where $d\lambda$ is the Liouville measure of $Z_1$, $\Phi_t = \phi_t\times \phi_t$ is the conjugacy of geodesic spaces, $\gamma'$ is the geodesic on $Y_t$ conjugate to $\gamma$ on $Z_1$, and $\theta'$ is the angle coordinate on $\mg(\gamma')$. Since $Z_1$ and $Z_2$ have the same length spectrum, $\Phi_t^* d\lambda$ is the Liouville measure of $Z_2$, and by prop \ref{spectraconverge}, the Liouville measures of $Y_t$ converge to this measure as $t\rightarrow 0$. Since the integral is compactly supported, the weak* uniform topology on $\mc(S)$ implies that there is some $\eta(t,a)$ so that
$$\frac{1}{l_{Z_1}(\gamma)}\int_\gamma \int_0^\pi F_a(\theta_t'(v,\theta))\sin\theta d\theta dt < \frac{l_{Y_t}(\gamma')}{l_{Z_1}(\gamma)}\int_0^\pi F_a(\theta')\sin\theta' d\theta' + \eta(t,a),$$
where $\eta(t,a) \rightarrow 0$ as $t\rightarrow 0$ or as $a\rightarrow 0$.

Now note that by the choice of the deformation, and since $Z_1$ and $Z_2$ have the same spectrum, there is some $M_t$ so that $\frac{l_{Y_t}(\gamma')}{l_{Z_1}(\gamma)} \leq M_t$ for any closed geodesic $\gamma$, and $M_t\rightarrow 1$ from above as $t\rightarrow 0$. That is, since the metrics on $Y_t$ are obtained by scaling the metric on $Z_2$ by an appropriate function, there is a limit to how much longer a geodesic on $Y_t$ can be than the isotopic geodesic on $Z_2$. Together with the approximation of $dv$ above by measures supported on closed geodesics, this shows that for any $a$,
$$\int_0^\pi F_a(\Theta_t'(\theta))\sin\theta d\theta \leq M_t \int_0^\pi F_a(\theta)\sin\theta d\theta + \eta(t,a),$$
where $M_t\rightarrow 1$ and $\eta(t,a)\rightarrow 0$ as $t \rightarrow 0$. Compare this to the similar hypothesis in lemma \ref{otallemma}. It remains to prove an adapted version of this lemma to show that $\Theta_t'$ is controllably close to the identity. For ease of notation, let $\psi = \Theta_t'$.

\vspace{10mm}
\begin{figure}[htb]
  \centering
  \def\svgwidth{250pt}
  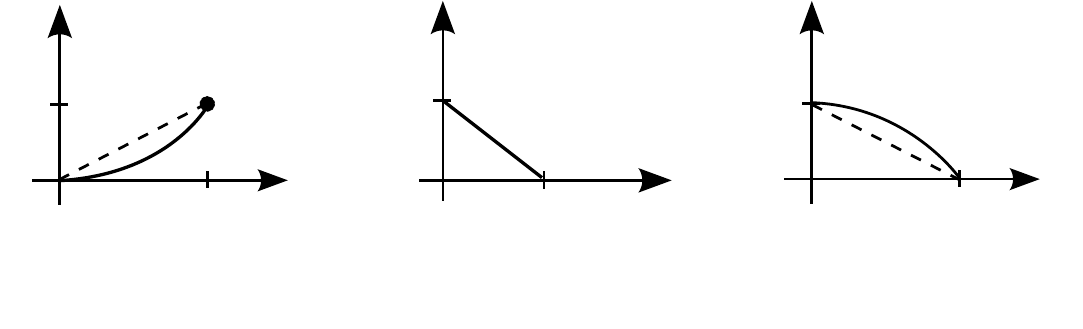
\end{figure}

Let $x$ be in $[0,\pi]$, with $\psi(x) = a$. Assume that $\psi(\theta) \leq \frac{a}{x}\theta$ for all $\theta \in [0,x]$. Then $F_a(\psi(\theta)) > -\frac{a}{x}\theta + a$, so:
\begin{eqnarray*}
a - a \frac{\sin x}{x} = \int_0^x \left(-\frac{a}{x}\theta + a\right)\sin\theta d\theta  &<& \int_0^\pi F_a(\psi(\theta))\sin\theta d\theta \\ \\
&\leq& M_t\int_0^\pi F_a(\theta)\sin\theta d\theta + \eta = M_t(a-\sin a) + \eta.
\end{eqnarray*}
Rearranging the outer two expressions of this inequality gives:
$$M_t \frac{\sin a}{a} - \left(M_t - 1 + \frac{\eta}{a}\right) < \frac{\sin x}{x}.$$

If $t\rightarrow 0$, then $M_t \rightarrow 1$ and $\eta\rightarrow 0$, so this inequality goes to $\frac{\sin a}{a} < \frac{\sin x}{x}$, which is equivalent to $a > x$, since $\frac{\sin x}{x}$ is decreasing on $[0,\pi]$. This corresponds to the statement in the proof of \ref{otallemma} that, when $\Theta'$ is comparing two structures with the same spectrum, there can be no initial interval on which $\Theta'$ is less than the identity.

For small $t$, this means that $a$ can be less than $x$, but there is a control on how much less. In other words, for each $x$, there is a slope $m_x$ so that $\psi(\theta)$ is not less than $m_x\theta$ on the entire interval $[0,x]$, and further each $m_x \rightarrow 1$ as $t\rightarrow 0$. 

Now fix a small $t$. As $x\rightarrow 0$, the $\eta$ term in the above inequality also goes to $0$, which implies that $m_x$ is increasing, and also bounded above by $1$. Let $m_t = \lim_{x\rightarrow 0} m_x$.

Then there is no interval $(0,x)$ so that $\psi(\theta) < m_t\theta$ on all of $(0,x)$. Indeed, if there were such an interval, then $\psi(\theta)$ would be less than one of the slopes that converge to $m_t$, which is impossible by the previous discussion. Thus there is a sequence of points $\{x_i\}$ converging to $0$ so that $\psi(x_i) \geq m_t x_i$ for all $i$.

\vspace{10mm}
\begin{figure}[htb]
  \centering
  \def\svgwidth{200pt}
  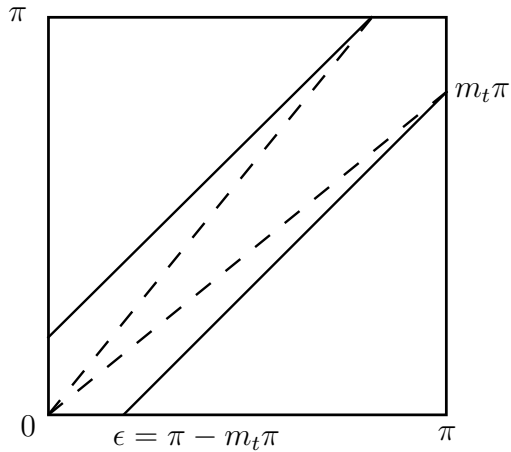
  \caption{Choosing $\epsilon$ so that $\psi$ is within $\epsilon$ of the identity.}
\end{figure}

Now let $\epsilon = \pi - m_t\pi$. Assume by absurd that there is some $b$ so that $|\psi(b) - b| > \epsilon$. By the symmetry of $\psi$ in $\pi - \theta$, we can assume that $\psi(b) < b - \epsilon$. Let $a = \sup\{\theta < b: \psi(\theta) > m_t\theta\}$. By the existence of the sequence $\{x_i\}$, this $\sup$ is not over an empty set. As in the proof of lemma \ref{otallemma}, $a < b$ and $\psi(a) = m_ta$. Choose some $x_i$ so that $a + x_i < b$. Then by the super-additivity of $\psi$,
$$m_ta + \psi(x_i) = \psi(a) + \psi(x_i) \leq \psi(a + x_i) < m_t(a + x_i) = m_ta + m_t x_i$$
This implies that $\psi(x_i) < m_t x_i$, which is a contradiction. Thus there is no such $b$, and $\psi$ is within $\epsilon$ of the identity on all of $[0,\pi]$. Since $m_t \rightarrow 1$ as $t\rightarrow 0$, we have also that $\epsilon \rightarrow 0$ as $t\rightarrow 0$. \qed

Let $\alpha, \beta, \gamma$ be three distinct non-singular geodesics in $\wt{Z_1}$ which share a common point, and let $T_t(\alpha,\beta,\gamma)$ be the triangle in $\wt{Y_t}$ formed by $\Phi_t(\alpha), \Phi_t(\beta), \Phi_t(\gamma)$. The previous proposition, along with the super-additivity argument on page \pageref{superadditive}, imply that the sum of the interior angles of $T_t(\alpha,\beta,\gamma)$ goes to $\pi$ as $t\rightarrow 0$. However, this is not enough to show that the area of the triangles goes to 0, since the curvature of the surface away from the cone points is also going to 0. By the Gauss-Bonnet formula for triangles, one would need to show that the sum $\Sigma_t$ of the interior angles goes to $\pi$ faster than $t$ goes to 0, so that $\frac{\pi-\Sigma_t}{t} \rightarrow 0$ as $t\rightarrow 0$. 
\appendix

\chapter{Further questions}

The space $\mc(S)$ of topological geodesic currents on $S$ is seen to be a useful ambient environment for studying moduli spaces of geometric structures on $S$. Otal's result on the spectral rigidity of Neg($S$) can be seen as saying that the mapping $L:\textrm{Neg}(S)\rightarrow \mc(S)$ is injective, so there is an embedding of Neg($S$) into $\mc(S)$. In fact, of course, there are many such embeddings possible, so it is important to realize what makes the Liouville map somehow the ``right" one. This is because of the property that $i(L_X, \alpha) = l_X(\alpha)$, for any $\alpha \in C$. 

In other words, the Liouville map is chosen so that the length spectrum of any Riemannian surface is recoverable from its current via the intersection form. Since Otal proved that currents are separated by their intersections with $C$ (see prop. \ref{otalcurrents}), the embedding can be seen as being uniquely determined by this property. We say that the Liouville map is then {\em length preserving.}

It is shown by Bonahon in \cite{Bonahon2} that there is also a length-preserving embedding of the cone of measured foliations $\mathcal{MF}(S)$ into $\mc(S)$ as the light cone of the intersection form, where the length of a foliation with respect to a metric is in the sense of Thurston.

It is conjectured herein that there is a length-preserving embedding of all of NonPos($S$) into $\mc(S)$, with the only possibility yet to be ruled out being that two different flat structures may define the same Liouville current. With this viewpoint in mind, we ask the following questions.

\section{Which functionals are currents?}

As in prop \ref{otalcurrents}, there is an embedding $I:\mc(S)\rightarrow \rc$ given by $\mu \mapsto (\alpha \mapsto i(\mu, \alpha))$. What is the image of this embedding? In other words, are there some algebraic or combinatorial properties of a functional $f:C\rightarrow\mathbf{R}$ that imply that $f = I(\mu)$ for some current $\mu$?

A partial answer is that such a functional $f = I(\mu)$ must be determined by its values on the primitive classes in $C$. That is, if $\alpha^n \in C$ is the class that represents a curve obtained by following a curve in $\alpha$ for $n$ periods, then necessarily $f(\alpha^n) = i(\mu, \alpha^n) = n i(\mu,\alpha) = n f(\alpha)$. Additionally, of course $f$ must take only non-negative values.

So if $C' \subset C$ is the subset of primitive curve classes, which functionals $C' \rightarrow \mathbf{R}_{\geq 0}$ extend to functionals in the image of $I$? It is clear that any one value can be chosen arbitrarily, and it may be possible to prove an inductive step that will allow any finite number of choices to be made. However, it seems unlikely that infinitely many values can be arbitrarily chosen for such a functional.

One reason this question is interesting is because there are types of structures on $S$ which define functionals in $\rc$ in some interesting way (see below on Hitchin components and cross-ratios), and it could be useful to know whether those functionals can be ``pulled back" to geodesic currents. There are few techniques for studying an arbitrary real-valued functional, but the theory of measures is rich. For instance, in Theorem 1 of this paper, two length spectra are shown to be distinct not by studying properties of the functionals themselves, but rather measure properties of their associated geodesic currents.

\section{Do length spectra converge between moduli spaces?}

On page \pageref{flatlimit} it is shown that a flat surface can be thought of as a limit of negatively curved cone surfaces. One very particular method of constructing such a limit is given, but there are many others. For instance, allowing the conformal structure or the area to vary would produce different types of deformations which limit to the same flat surface. In section \ref{flatrigidity} it is shown that for this particular well-chosen construction, the length spectrum along the deformation converges to the spectrum of the flat limit surface. Does this happen regardless of the deformation chosen?

Furthermore, a negatively-curved cone surface can be thought of as a limit of negatively curved Riemannian surfaces. To introduce a cone point into a Riemannian structure $X$, take a metric ball $B_\epsilon$ centered at the point, and let $k$ be the total curvature inside $B_\epsilon$ (note that $k<0$). Letting $\epsilon \rightarrow 0$, while keeping the total curvature in the ball constant and leaving the metrics on each $X \setminus B_\epsilon$ isometric, produces in the limit a point with cone angle $2\pi - k$ (i.e. with concentrated curvature $k$).

Similarly, a flat surface can be a limit of negatively curved Riemannian surfaces, if all the curvature of the surface is gradually moved into small neighborhoods of the cone points as the size of these neighborhoods goes to 0.

Do all of these types of limits also cause the length spectra to converge? If this were true, then by the injectivity of the map $I: \mc(S)\rightarrow \rc$ and the completeness of $\mc(S)$, the Liouville currents would also converge. This would mean that Neg*($S$) sits in $\mc(S)$ as a sort of partial boundary to Neg($S$), and the image of Flat*($S$) would similarly be a partial boundary for both Neg($S$) and Neg*($S$). It would be interesting to know how these moduli spaces all fit together as currents. 

\section{Are Hitchin components representable as geodesic currents?}

For any $n \geq 2$, a representation $\Gamma \rightarrow PSL(n, \mathbf{R})$ is called n-Fuchsian if it can be written as a composition of a Fuchsian representation $\Gamma\rightarrow PSL(2,\mathbf{R})$ with the irreducible representation $PSL(2,\mathbf{R})\rightarrow PSL(n,\mathbf{R})$. A representation is called n-Hitchin if it can be deformed into an n-Fuchsian representation. Hitchin proved in \cite{Hitchin} that for odd $n$ there is a single component of such representations, and for even $n$ two isomorphic ones, each homeomorphic to a ball.

In \cite{Labourie1}, Labourie shows that any Hitchin representation $\rho$ is into matrices of split real type, and uses this to define a {\em period} functional $\omega_\rho: \Gamma \rightarrow \mathbf{R}$, where $\omega_\rho(\gamma)$ is the log of the ratio of the largest and smallest eigenvalues of $\rho(\gamma)$. Since this is a conjugacy invariant, the period can be thought of as a mapping from each Hitchin component $H_n \rightarrow \rc$. For any $n$, is this mapping into the image of $\mc(S)$? In other words, are Hitchin components naturally representable as geodesic currents, in such a way that the intersection form produces the periods of the representations?

For $n=2$, this is trivially true, since the two 2-Hitchin components are simply the Fricke space of holonomies of hyperbolic structures. For $n=3$, Choi and Goldman showed in \cite{ChoiGoldman} that the 3-Hitchin component consists of convex real projective structures on $S$, with the periods corresponding to the Hilbert length spectrum. For the Fuchsian representations, these are equivalent to hyperbolic structures with the hyperbolic length spectrum, so these representations are clearly representable as currents (and overlap the 2-Hitchin component). Will this extend to the quasi-Fuchsian 3-Hitchin representations?

Labourie conjectures in \cite{Labourie2} that the union of all the images of the Hitchin components in $\rc$ contains the image of Neg($S$). It would be interesting to know how all of these spaces overlap and interact.

\renewcommand{\baselinestretch}{1}
\small\normalsize

\newpage
\bibliographystyle{plain}
\bibliography{thesisbib}

\end{document}